\documentclass[12pt]{article}
\usepackage{amssymb}
\usepackage{latexsym}

\newtheorem{definition}{Definition}[section]
\newtheorem{theorem}[definition]{Theorem}
\newtheorem{lemma}[definition]{Lemma}
\newtheorem{corollary}[definition]{Corollary}

\newtheorem{example}[definition]{Example}

\newtheorem{problem}[definition]{Problem}
\newtheorem{note}[definition]{Note}

\newtheorem{proposition}[definition]{Proposition}
\newtheorem{notation}[definition]{Notation}

\def\I{\mathbb I}

\def\R{\mathbb R}
\def\C{\mathbb C}

\def\Z{\mathbb Z}
\def\F{\mathbb F}
\def\K{\mathbb F}

\begin{document}

\title{\bf Finite-dimensional irreducible modules  \\
for the three-point  
$\mathfrak{sl}_2$ loop algebra
}
\author{
Tatsuro Ito{\footnote{
Department of Computational Science,
Faculty of Science,
Kanazawa University,
Kakuma-machi,
Kanazawa 920-1192, Japan
}}
{\footnote{Supported in part by JSPS grant 18340022}}
$\;$ and
Paul Terwilliger{\footnote{
Department of Mathematics, University of
Wisconsin, 480 Lincoln Drive, Madison WI 53706-1388 USA}
}}
\date{}
%to get date printout, comment out above line

\maketitle
\begin{abstract}
Recently Brian Hartwig and the second author found a presentation
for the three-point 
$\mathfrak{sl}_2$ loop algebra by generators and relations.
To obtain this presentation they defined a Lie algebra
$\boxtimes$ by generators and relations, and displayed an
isomorphism  from
$\boxtimes$ to the
three-point 
$\mathfrak{sl}_2$ loop algebra. In this paper we 
describe the finite-dimensional irreducible
$\boxtimes$-modules from multiple points of view.

\bigskip
\noindent
{\bf Keywords}. 
Tetrahedron Lie algebra,
tridiagonal pair,
 Onsager Lie algebra,
Kac-Moody algebra.
 \hfil\break
\noindent {\bf 2000 Mathematics Subject Classification}. 
Primary: 17B37. Secondary: 16W35,
05E35, 
82B23.
 \end{abstract}

\section{Introduction}
In \cite{HT}
Hartwig and the second author
found a presentation for the three-point $\mathfrak{sl}_2$ loop
algebra via generators and relations. To obtain this
presentation they defined a Lie algebra $\boxtimes$ 
(pronounced  ``tet'')
by
generators and relations, and displayed an isomorphism
 from $\boxtimes$ to
 the three-point $\mathfrak{sl}_2$ loop algebra.
 $\boxtimes$ has essentially six generators, and it is natural
 to identify these with the six edges of a tetrahedron.
The action of the symmetric group $S_4$
on the tetrahedron
 induces an action of $S_4$ on
$\boxtimes$ as a group of automorphisms \cite[Section 2]{HT}.
For each face of the
 tetrahedron the three
 surrounding edges
 form a basis for a subalgebra of $\boxtimes$
 that is isomorphic to $\mathfrak{sl}_2$ \cite[Corollary 12.4]{HT}.
 Any five of the six edges of the tetrahedron generate a subalgebra of
 $\boxtimes$
 that is isomorphic to the $\mathfrak{sl}_2$ loop algebra
  \cite[Corollary 12.6]{HT}.
  Each pair of opposite
  edges of the tetrahedron
  generate a subalgebra of $\boxtimes$ that is isomorphic
  to the Onsager algebra $\mathcal O$ \cite[Corollary 12.5]{HT};
  let us call these Onsager subalgebras.
  Then $\boxtimes $
  is the direct sum of its three Onsager subalgebras
   \cite[Theorem 11.6]{HT}.
   In \cite{E} Elduque found an
   attractive decomposition of $\boxtimes$ into a direct sum of three
   abelian subalgebras, and he
   showed how these subalgebras
   are related to the Onsager subalgebras.
   In \cite{PT} Pascasio and the second author gave an
   action of $\boxtimes$ on the standard module for each
   Hamming association
   scheme. 
  % In
  % \cite{ITdrg} we described a 
  % $q$-analog of this phenomenon.
   In \cite{Br} Bremner
   obtained
   the universal central extension of
   the three-point $\mathfrak{sl}_2$ loop algebra.
   By modifying the defining relations for $\boxtimes$,
   Benkart and the second author obtained
   a presentation for this
   extension by generators and relations
    \cite{BT}.
    In \cite{Ha} Hartwig classified the
    finite-dimensional irreducible $\boxtimes$-modules
    over an algebraically closed field $\F$ with characteristic 0.
    He did this by displaying a bijection between
    (i) the set of isomorphism classes of finite-dimensional 
    irreducible $\boxtimes$-modules; (ii) the set of
    isomorphism classes of finite-dimensional irreducible
    $\mathcal O$-modules that have type $(0,0)$
   \cite[Theorems~1.7,~1.8]{Ha}. 
   The modules in (ii) were classified
   earlier by Davies
\cite{Davfirst},
\cite{Da}; see also
Date and Roan
\cite{DateRoan2}. A summary of this classification is
given in
   \cite[Theorems~1.3,~1.4,~1.6]{Ha}.

\medskip
\noindent 
In this paper we consider further the 
finite-dimensional irreducible $\boxtimes$-modules.
We pick up where
Davies, Date and Roan, and Hartwig left off,
and obtain a description of these modules
that has more detail and takes a different point of view.
Our results are summarized as follows.
Among the finite-dimensional irreducible
$\boxtimes$-modules we identify a special case 
called an evaluation module. 
Each evaluation module is written
$V_d(a)$, where
$d\geq 1$ is an integer and
 $a\in \F\backslash\lbrace 0,1\rbrace$.
The scalar $a$ is called the evaluation parameter.
To get the $\boxtimes$-module $V_d(a)$
we start with the irreducible 
$\mathfrak{sl}_2$-module $V_d$ of 
dimension $d+1$, and pull back the
$\mathfrak{sl}_2$-module structure
via an
 evaluation homomorphism
$EV_a :\boxtimes \to 
\mathfrak{sl}_2$.
We note that our evaluation homomorphism is different
from the one in
\cite{DateRoan2}.   
We show that every finite-dimensional irreducible
$\boxtimes$-module is a tensor product of 
evaluation modules.
We give three characterizations that
enable us to
recognize the evaluation modules  
among all finite-dimensional irreducible
$\boxtimes$-modules.
Upon twisting an evaluation module via an element
of $S_4$ we get an evaluation module whose
evaluation parameter is potentially different.
We describe how twisting
affects the evaluation parameter.
To do this, 
we first display an action of $S_4$
on 
$\F\backslash\lbrace 0,1\rbrace$ by linear
fractional transformations.
We then show that for an evaluation module
$V_d(a)$ and 
$\sigma \in S_4$, the $\boxtimes$-module
$V_d(a)$ twisted via
$\sigma $ is isomorphic to
$V_d(\sigma(a))$.
Let $G$ denote the kernel of the above
$S_4$-action on
$\F\backslash\lbrace 0,1\rbrace$.
We show that $G$ is the unique normal subgroup of
$S_4$ that has cardinality 4.
For a given evaluation module $V_d(a)$ we consider 
24 bases that are described as follows.
Let $\I$ denote the vertex set of the tetrahedron.
For mutually distinct $i,j,k,\ell \in \I$ we
define an (ordered) basis
$\lbrack i,j,k,\ell\rbrack$ of
 $V_d(a)$ such that
(i) the basis diagonalizes the $\boxtimes$-generator
identified with the edge of the tetrahedron incident with $k$ and $\ell$;
(ii) the sum of the basis vectors is a common eigenvector
for the three $\boxtimes$-generators identified with the edges 
 incident with $i$. 
We find the matrices that represent all the 
$\boxtimes$-generators with respect to
$\lbrack i,j,k,\ell\rbrack$.
Moreover we find
the transition matrices from the basis
$\lbrack i,j,k,\ell\rbrack $ to the bases
\begin{eqnarray*}
\lbrack j,i,k,\ell\rbrack,
\qquad \qquad 
\lbrack i,k,j,\ell\rbrack,
\qquad \qquad 
\lbrack i,j,\ell,k\rbrack.
\end{eqnarray*}
The first matrix is diagonal, 
the second matrix is lower triangular, and the third
matrix is the identity reflected about a vertical axis.
We wind up our discussion of
evaluation modules by showing how
they can be concretely realized using
homogeneous polynomials in two variables.
Now consider a general
finite-dimensional irreducible
$\boxtimes$-module $V$.
 For $\sigma \in G$ we show that
$V$ twisted via $\sigma$ is isomorphic to
$V$. We associate with $V$ two polynomials whose
algebraic structure reflects that of $V$. 
The first polynomial, denoted $S_V$ and
called the shape polynonial, has coefficients
that encode
the eigenspace dimensions for the action of the
$\boxtimes$-generators on $V$.
We show how a certain factorization of $S_V$  
over the integers
corresponds to the factorization of $V$ into a tensor
product of evaluation modules.
The second polynomial, denoted $P_V$ and
called the Drinfel'd polynomial, is analogous
to the Drinfel'd polynomial for the classical $\mathfrak{sl}_2$
loop algebra 
\cite{charint},
\cite{charp},
\cite{degxz}.
We show that the factorization of $P_V$ 
into linear factors corresponds to the
factorization of $V$ into a tensor product of
evaluation modules.
Moreover we show that the map $V\mapsto P_V$ induces a bijection
between (i) the set of isomorphism classes of
finite-dimensional irreducible $\boxtimes$-modules;
(ii) the set of univariable polynomials over $\F$
that have constant coefficient 1 and do not vanish at 1.
We display 
a nondegenerate bilinear form $\langle\,,\,\rangle$ on $V$ such that
\begin{eqnarray*}
\langle \xi.u,v\rangle = -\langle u,\xi.v\rangle
\qquad \qquad \xi \in \boxtimes, \qquad u,v\in V.
\end{eqnarray*}
We use this form to show that the
$\boxtimes$-module $V$ is isomorphic to
the dual $\boxtimes$-module $V^*$.
There is another type of bilinear form on $V$ that
is of interest.
For a nonidentity $\sigma \in G$ we 
display a nondegenerate symmetric bilinear
form
$\langle\,,\,\rangle_\sigma$ on $V$ such that
\begin{eqnarray*}
\langle \xi.u,v\rangle_\sigma = 
- 
\langle u,\sigma(\xi).v\rangle_\sigma
\qquad \qquad 
\xi \in \boxtimes, \qquad u,v \in V.
\end{eqnarray*}
We remark on our motivation to
investigate
$\boxtimes$. This has to do with
the tridiagonal pairs
\cite{TD00},
\cite{shape},
\cite{tdanduq},
\cite{NN}
and the closely related Leonard pairs
\cite{nom2},
\cite{nom5},
\cite{nom6},
\cite{LS99},
  \cite{qSerre},
    \cite{lsint},
\cite{qrac},
\cite{aw}.
A Leonard pair is a pair of semisimple
linear transformations on a finite-dimensional
vector space, each of which acts
in an irreducible tridiagonal fashion
on an eigenbasis for the other
\cite[Definition 1.1]{LS99}.
There is a close connection between the
Leonard pairs and the orthogonal polynomials
that make up the terminating branch of the
Askey scheme 
\cite{KoeSwa},
\cite{LS99},
\cite{TLT:array}.
A tridiagonal pair is a mild generalization
of a Leonard pair 
\cite[Definition 1.1]{TD00}.
A pair of 
$\boxtimes$-generators identified with
opposite edges of the tetrahedron
act 
on each finite-dimensional irreducible
$\boxtimes$-module as a tridiagonal pair 
 \cite[Theorem 1.7,~Corollary 2.7]{Ha}.
For the details on this connection
see 
 \cite[Section~2]{Ha}.

\section{The three-point 
$\mathfrak{sl}_2$ loop algebra} 

\noindent
In this section we recall the 
three-point 
$\mathfrak{sl}_2$ loop algebra and its
relationship to $\boxtimes$.

\medskip
\noindent Throughout the paper $\K$ denotes an algebraically
closed field with characteristic 0.
By $\mathfrak{sl}_2$ we mean the Lie algebra
$\mathfrak{sl}_2(\K)$.

\begin{definition}
\label{def:3point}
\rm
Let $t$ denote an indeterminate and let
$\F\lbrack t,t^{-1}, (t-1)^{-1}\rbrack$ denote the
$\F$-algebra of all Laurent polynomials in $t,t-1$ that
have coefficients in $\F$. We abbreviate
\begin{eqnarray*}
\mathcal{A}=\F\lbrack t,t^{-1}, (t-1)^{-1}\rbrack.
\end{eqnarray*}
We consider the
Lie algebra over $\K$ consisting of the $\F$-vector space
\begin{eqnarray*}
 {\mathfrak{sl}}_2 \otimes \mathcal{A}, 
\qquad \qquad \otimes=\otimes_{\K}
\end{eqnarray*}
 and Lie bracket
\begin{eqnarray*}
\lbrack u\otimes a,v\otimes b\rbrack =
\lbrack u,v\rbrack \otimes ab,
\qquad \qquad u,v \in
 {\mathfrak{sl}}_2,
\qquad a,b \in 
\mathcal{A}.
\end{eqnarray*}
We call 
 ${\mathfrak{sl}}_2 \otimes \mathcal{A}$ the
 {\it three-point
 ${\mathfrak{sl}}_2$  loop algebra}.
\end{definition}
\noindent See \cite{Br}, \cite{sch}, \cite{sch2}, \cite{sch3}
for information on multipoint loop algebras and
related topics.

\medskip
\noindent We now recall the Lie
algebra  $\boxtimes$.

\begin{definition}
\label{def:tet}
\rm
\cite[Definition 1.1]{HT}
Let $\boxtimes$ denote the Lie algebra over $\K$
that has  generators
\begin{eqnarray}
\label{eq:boxgen}
\lbrace x_{ij} \,|\,i,j\in \I, i\not=j\rbrace
\qquad \qquad \I = \lbrace 0,1,2,3\rbrace
\end{eqnarray}
and the following relations:
\begin{enumerate}
\item[{\rm (i)}]  
For distinct $i,j\in \I$,
\begin{eqnarray*}
x_{ij}+x_{ji} = 0.
\label{eq:rel0}
\end{eqnarray*}
\item[{\rm (ii)}]  
For mutually distinct $i,j,k\in \I$,
\begin{eqnarray*}
\lbrack x_{ij},x_{jk}\rbrack = 2x_{ij}+2x_{jk}.
\label{eq:rel1}
\end{eqnarray*}
\item[{\rm (iii)}]  
For mutually distinct $i,j,k,\ell \in \I$,
\begin{eqnarray}
\label{eq:dg}
\lbrack x_{ij},
\lbrack x_{ij},
\lbrack x_{ij},
x_{k \ell}\rbrack \rbrack \rbrack= 
4 \lbrack x_{ij},
x_{k \ell}\rbrack.
\label{eq:rel2}
\end{eqnarray}
\end{enumerate}
We call $\boxtimes$ the {\it tetrahedron algebra}
or ``tet'' for short. 
\end{definition}

\begin{note} \rm
The relations (\ref{eq:dg})
are known as the {\it Dolan-Grady} relations
\cite{CKOnsn},
\cite{DateRoan2},
\cite{Davfirst},
\cite{Da},
\cite{Dolgra},
\cite{Klish3}.
\end{note}

\noindent
After a comment we will display a Lie algebra
isomorphism from $\boxtimes$ to
$\mathfrak{sl}_2 \otimes \mathcal{A}$.

\medskip
\noindent
 Recall that 
the Lie algebra $\mathfrak{sl}_2$
has a basis $e,f,h$ 
 and
\begin{eqnarray*}
\lbrack h,e\rbrack = 2e,
\qquad \qquad
\lbrack h,f\rbrack = -2f,
\qquad
\qquad
\lbrack e,f\rbrack = h.
\end{eqnarray*}
Following \cite[Lemma 3.2]{HT} we define
\begin{eqnarray}
x=2e-h,
\qquad \qquad 
y=-2f-h,
\qquad \qquad 
z=h.
\label{eq:xyz}
\end{eqnarray}
Then $x,y,z$ is a basis for 
$\mathfrak{sl}_2$ and 
\begin{eqnarray*}
\label{eq:equit}
\lbrack
x,y
\rbrack =  2x+2y,
\qquad 
\qquad
\lbrack
y,z
\rbrack =  2y+2z,
\qquad 
\qquad
\lbrack
z,x
\rbrack =  2z+2x.
\end{eqnarray*}
By analogy with \cite{equit1},
\cite{equit2}
we call 
$x,y,z$ the {\it equitable basis} for 
$\mathfrak{sl}_2$. See \cite{equitsl2} for
a detailed study of this basis.

\begin{proposition}
\label{prop:eta}
{\rm \cite[Proposition 6.5]{HT}}
There exists an isomorphism of Lie algebras
$\psi: \boxtimes \to \mathfrak{sl}_2 \otimes \mathcal{A}$
that sends
\begin{eqnarray*}
&&x_{12}\mapsto x\otimes 1, \qquad \qquad x_{03} \mapsto y\otimes t + z\otimes (t-1),
\\
&&x_{23}\mapsto y\otimes 1, \qquad \qquad x_{01} \mapsto z\otimes(1-t^{-1})
- x\otimes t^{-1},
\\
&&x_{31}\mapsto z\otimes 1, \qquad \qquad x_{02} \mapsto x\otimes(1-t)^{-1}
+ y\otimes t(1-t)^{-1},
\end{eqnarray*}
where $x,y,z$ is the equitable basis for 
$\mathfrak{sl}_2$.
\end{proposition}

\noindent We mention a result for later use.

\begin{lemma} {\rm \cite[Corollary 12.1]{HT}}
\label{lem:sl2inj}
For mutually distinct $i,j,k \in \I$
there exists an injective homomorphism
of Lie algebras from
$\mathfrak{sl}_2$ to $\boxtimes$ that sends
\begin{eqnarray*}
x \;\mapsto \; x_{ij},
\qquad \qquad 
y \;\mapsto  \;x_{jk},
\qquad \qquad 
z \;\mapsto  \;x_{ki}.
\end{eqnarray*}
\end{lemma}

\section{Finite-dimensional irreducible $\boxtimes$-modules}

\noindent In this section we review some basic
facts and notation concerning
finite-dimensional irreducible $\boxtimes$-modules.
This material is summarized from
\cite{Ha}.

\medskip
\noindent 
Let $V$ denote a vector space over $\K$ with finite positive
dimension.
Let $\lbrace s_n\rbrace_{n=0}^d $ denote a sequence
of positive integers whose sum is the dimension of $V$.
By a {\it decomposition of $V$ of shape 
$\lbrace s_n\rbrace_{n=0}^d$}
we mean a sequence 
$\lbrace V_n\rbrace_{n=0}^d$ of subspaces of $V$ such that
$V_n$ has dimension $s_n$ for $0 \leq n \leq d$ and $V=\sum_{n=0}^d V_n$
(direct sum). We call $V_n$ the $n$th {\it component} of
the decomposition. We call $d$ the {\it diameter} of the
decomposition.
For notational convenience we define
$V_{-1}=0$ and $V_{d+1}=0$. By the {\it inversion} of
$\lbrace V_n\rbrace_{n=0}^d$ we mean the decomposition
$\lbrace V_{d-n}\rbrace_{n=0}^d$.

\medskip
\noindent 
Now let $V$ denote a finite-dimensional irreducible
$\boxtimes$-module. By 
\cite[Theorem 3.8]{Ha}
each generator $x_{ij}$ of $\boxtimes$ is semisimple
on $V$. Moreover 
there exists an integer $d\geq 0$ such that
for each generator $x_{ij}$ the set of distinct eigenvalues
on $V$ is $\lbrace d-2n\,|\,0 \leq n \leq d\rbrace$. 
We call
$d$ the {\it diameter} of $V$.
For distinct $i,j \in \I$ we define a decomposition of $V$
called $\lbrack i,j\rbrack$. The decomposition
$\lbrack i,j\rbrack$ has diameter $d$.
For $0 \leq n \leq d$ the $n$th component 
of $\lbrack i,j\rbrack$ 
is the eigenspace of $x_{ij}$ for the eigenvalue $2n-d$.
Using Definition
\ref{def:tet}(i) we find $\lbrack j,i\rbrack$ is
the inversion of $\lbrack i,j\rbrack$.
By \cite[Corollary 3.6]{Ha}, for distinct $i,j \in \I$
the shape of 
$\lbrack i,j\rbrack$ is independent of $i,j$. We denote this
shape by $\lbrace \rho_n \rbrace_{n=0}^d$ and note that
$\rho_n=\rho_{d-n}$ for $0 \leq n\leq d$.
By the {\it shape of $V$} we mean the sequence
$\lbrace \rho_n \rbrace_{n=0}^d$.

\medskip
\noindent 
Let $V$ denote a finite-dimensional irreducible $\boxtimes$-module.
For distinct $i,j \in \I$ and for distinct $r,s \in \I$ 
we now describe the action of 
$x_{rs}$ on the decomposition 
$\lbrack i,j\rbrack$ of $V$. 
Denote this decomposition by 
$\lbrace V_n\rbrace_{n=0}^d$. Then by
\cite[Theorem 3.3]{Ha},
for $0 \leq n \leq d$
the action of $x_{rs}$ on $V_n$ is given in the table below.

\bigskip

\centerline{
\begin{tabular}[t]{c|c}
       {\rm Case} & {\rm Action of $x_{rs}$ on $V_n$}
 \\ \hline  \hline
	$r=i, \;\quad s=j$  & $(x_{rs}-(2n-d)I)V_n=0$    
	\\
	$r=j, \;\quad s=i$  & $(x_{rs}-(d-2n)I)V_n=0$    
        \\
	$r=j, \;\quad s\not=i$  & $(x_{rs}-(d-2n)I)V_n\subseteq V_{n+1}$ 
	\\
	$r\not=i,\; \quad s=j$  & $(x_{rs}-(2n-d)I)V_n\subseteq V_{n+1}$ 
	\\
	$r=i, \;\quad s\not=j$  & $(x_{rs}-(2n-d)I)V_n\subseteq V_{n-1}$ 
	\\
	$r\not=j,\; \quad s=i$  & $(x_{rs}-(d-2n)I)V_n\subseteq V_{n-1}$ 
	\\
	$i,j,r,s$ distinct & $x_{rs}V_n \subseteq V_{n-1}+V_n+V_{n+1}$
\end{tabular}}

\bigskip
\noindent We recall the notion of a {\it flag}.
For the moment let $V$ denote a vector space over $\K$ with finite
positive dimension and let $\lbrace s_n\rbrace_{n=0}^d$ denote
a  sequence of positive integers whose sum is the dimension
of $V$. By a {\it flag on $V$ of shape 
$\lbrace s_n\rbrace_{n=0}^d$} we mean
a nested sequence $U_0\subseteq U_1\subseteq \cdots \subseteq
U_d$ of subspaces of $V$ such that the dimension of $U_n$
is $s_0+\cdots+s_n$ for $0 \leq n \leq d$. We call
$U_n$ the $n$th {\it component} of the flag. We call
$d$ the {\it diameter} of the flag. We observe 
that $U_d=V$.
The following construction yields a flag on $V$.
Let $\lbrace V_n \rbrace_{n=0}^d$ denote a decomposition of $V$
of shape 
 $\lbrace s_n \rbrace_{n=0}^d$.
 Define
\begin{eqnarray*}
U_n=V_0+V_1+\cdots+V_n \qquad \qquad (0 \leq n\leq d).
\end{eqnarray*}
Then the sequence 
$U_0 \subseteq U_1 \subseteq \cdots \subseteq U_d$
is a flag on $V$ of shape
$\lbrace s_n \rbrace_{n=0}^d$.
We say this flag is {\it induced} by the decomposition
$\lbrace V_n \rbrace_{n=0}^d$.
We now recall what it means for two flags to be {\it opposite}.
Suppose we are given two flags on $V$ with the same diameter: 
 $U_0 \subseteq U_1 \subseteq \cdots \subseteq U_d$
and 
 $U'_0 \subseteq U'_1 \subseteq \cdots \subseteq U'_d$.
We say these flags are {\it opposite} whenever there exists
a decomposition 
$\lbrace V_n \rbrace_{n=0}^d$
 of $V$ such that
\begin{eqnarray*}
U_n=V_0+V_1+\cdots+V_n, \qquad \qquad 
U'_n=V_d+V_{d-1}+\cdots+V_{d-n}
\end{eqnarray*}
for $0 \leq n\leq d$. In this case
\begin{eqnarray}
\label{eq:zero}
U_i\cap U'_j = 0 \qquad \mbox{if}\quad i+j<d \qquad \qquad (0 \leq i,j\leq d)
\end{eqnarray}
and
\begin{eqnarray}
\label{eq:recover}
V_n=U_n\cap U'_{d-n} \qquad \qquad (0 \leq n\leq d).
\end{eqnarray}
In particular the decomposition
$\lbrace V_n \rbrace_{n=0}^d$
is uniquely
determined by the given flags.

\medskip
\noindent 
 We now return our attention to 
$\boxtimes$-modules.
Let $V$ denote a finite-dimensional irreducible $\boxtimes$-module
of diameter $d$. By \cite[Lemma 5.3]{Ha}
there exists a collection of
flags on $V$, denoted $\lbrack i \rbrack, i \in \I$, 
such that for distinct $i,j\in \I$ the
 decomposition
$\lbrack i,j\rbrack$ induces $\lbrack i\rbrack$.
By construction, for $i \in \I$ the 
 the shape
of the flag $\lbrack i \rbrack $ coincides with the shape of $V$.
By construction,
%\cite[Lemma 5.6]{Ha}
the flags $\lbrack i \rbrack$, $i \in \I$ are mutually
opposite.
Also by construction,
%\cite[Lemma 5.7]{Ha},
for distinct $i,j\in \I$ and for
$0 \leq n\leq d$ the $n$th component of $\lbrack i,j \rbrack$
is equal to the intersection of the following two sets:
\begin{enumerate}
\item[{\rm (i)}]  
component $n$ of the flag $\lbrack i\rbrack$;
\item[{\rm (ii)}] 
component $d-n$ of the flag $\lbrack j\rbrack$.
\end{enumerate}

\section{Evaluation modules for $\boxtimes$; preliminaries}

\noindent In the next section we will define
a class of $\boxtimes$-modules called evaluation modules.
In this section we make some preliminary comments
to set the stage.

\begin{lemma}
\label{lem:ideal}
Let $K$ denote an ideal of $\boxtimes$ such that $K\not=\boxtimes$.
Then for mutually distinct $i,j,k \in \I$ the following
sum is direct:
\begin{eqnarray}
\label{eq:ds}
\K x_{ij}+
\K x_{jk}+
\K x_{ki} + K.
\end{eqnarray}
\end{lemma}
\noindent {\it Proof:} 
We first claim that $x_{rs} \not\in K$ for all distinct 
$r,s \in \I$. To prove the claim
it suffices to show that $x_{rs}\in K$ implies $x_{su} \in K$
for  mutually distinct $r,s,u \in \I$.
If $x_{rs} \in K$ then 
$\lbrack x_{rs},x_{su}\rbrack \in K$ since $K$ is an ideal.
Also 
$\lbrack x_{rs},x_{su}\rbrack = 2x_{rs}+2x_{su}$
by 
 Definition
\ref{def:tet}(ii)
so $x_{su} \in K$ and
the claim follows. We can now easily show
that the sum 
(\ref{eq:ds}) is direct.
By Lemma
\ref{lem:sl2inj} the generators
$x_{ij}$, 
$x_{jk}$, 
$x_{ki}$ form a basis for
a subalgebra $L$ of $\boxtimes$ that is isomorphic to
$\mathfrak{sl}_2$.
Recall that 
$\mathfrak{sl}_2$ is simple so $L$ is simple.
By this and since $K\cap L$ is an ideal of $L$ we
find $K\cap L=0$ or $L \subseteq K$. The second possibility
cannot occur by the claim,
 so
$K\cap L=0$ and this means
that the sum (\ref{eq:ds}) is direct.
\hfill $\Box $ \\ 

\begin{lemma}
\label{lem:image}
For a nonzero Lie algebra homomorphism $\boxtimes \to \mathfrak{sl}_2$
and for mutually distinct $i,j,k \in \I$ the images of
$x_{ij}$, 
$x_{jk}$, 
$x_{ki}$  form a basis for
$\mathfrak{sl}_2$.
\end{lemma}
\noindent {\it Proof:} 
Let $K$ denote the kernel of the homomorphism and note
that $K$ is an ideal of $\boxtimes$. The homomorphism
is nonzero so $K\not=\boxtimes$. Applying Lemma
\ref{lem:ideal}
to $K$ 
we find that the images of
$x_{ij}$, 
$x_{jk}$, 
$x_{ki}$ are linearly independent in 
$\mathfrak{sl}_2$. 
 The result follows since
$\mathfrak{sl}_2$ has dimension 3.
\hfill $\Box $ \\ 

\begin{proposition}
\label{prop:a}
For a nonzero Lie algebra homomorphism
$\boxtimes \to \mathfrak{sl}_2$ there exists
a unique $a \in \K$ such that the following
are in the kernel:
\begin{eqnarray*}
&&
a x_{01}+(1-a)x_{02}-x_{03},
\qquad \qquad 
a x_{10}+(1-a)x_{13}-x_{12},
\\
&&
a x_{23}+(1-a)x_{20}-x_{21},
\qquad \qquad
a x_{32}+(1-a)x_{31}-x_{30}.
\end{eqnarray*}
Moreover $a\not=0$ and $a\not=1$.
\end{proposition}
\noindent {\it Proof:} 
We first claim that for mutually distinct 
$i,j,k,\ell \in \I$
the images of 
\begin{eqnarray*}
x_{ij}-x_{ik}, \qquad x_{ik}-x_{i\ell}
\end{eqnarray*}
in 
$\mathfrak{sl}_2$
are nonzero and linearly dependent.
These images are nonzero by Definition 
\ref{def:tet}(i) and
Lemma
\ref{lem:image}.
To show that they are linearly dependent,
let $K$ denote the kernel of the homomorphism.
By Lemma
\ref{lem:image} there exist scalars $\alpha, \beta, \gamma \in \K$
such that
\begin{eqnarray}
\label{eq:xijdep}
x_{i\ell}-\alpha x_{ij}-\beta x_{jk} - \gamma x_{ki}
\end{eqnarray}
is in $K$. Taking the Lie bracket of 
(\ref{eq:xijdep}) with $x_{ij}$ and simplifying the result
using Definition
\ref{def:tet}(i),(ii) we find
\begin{eqnarray}
\label{eq:xijdep2}
x_{i\ell} + (\beta-\gamma-1)x_{ij}+\beta x_{jk}-\gamma x_{ki}
\end{eqnarray}
is in $K$. Comparing  
(\ref{eq:xijdep}), 
(\ref{eq:xijdep2}) and using Lemma
\ref{lem:image} 
we find $\beta=0$ and $\gamma=\alpha-1$. Evaluating 
(\ref{eq:xijdep}) using this and Definition \ref{def:tet}(i)
we find that $K$ contains
\begin{eqnarray*}
\alpha (x_{ij}-x_{ik})+x_{ik}-x_{i\ell}
\end{eqnarray*}
and the claim follows.
We now apply the claim for the following
values of 
$(i,j,k,\ell)$:
\begin{eqnarray*}
(0,1,2,3),
\qquad \qquad
(1,0,3,2),
\qquad \qquad
(2,3,0,1),
\qquad \qquad
(3,2,1,0).
\end{eqnarray*}
We find that there exist unique
scalars $a,b,c,d \in \K$ such that each of
\begin{eqnarray}
\label{eq:K1}
&&
ax_{01}+(1-a)x_{02}-x_{03},
\\
&&
bx_{10}+(1-b)x_{13}-x_{12},
\label{eq:K2}
\\
&&
cx_{23}+(1-c)x_{20}-x_{21},
\label{eq:K3}
\\
&&
dx_{32}+(1-d)x_{31}-x_{30}
\label{eq:K4}
\end{eqnarray}
is in $K$.
None of $a,b,c,d$ is equal to 0 or 1 in
view of Lemma
\ref{lem:image} and Definition
\ref{def:tet}(i).
Define $E \in \boxtimes$ to be
(\ref{eq:K1}) plus $a/b$ times
(\ref{eq:K2})
plus $(1-a)/(1-c)$ times 
(\ref{eq:K3}) plus 
(\ref{eq:K4}).
By construction $E \in K$.
Also by Definition \ref{def:tet}(i),
$E$ 
is equal to
a certain linear combination
of $x_{12}, x_{23},x_{31}$. 
In this linear combination each coefficient
must be 0 since
the images of $x_{12}, x_{23}, x_{31}$
are linearly independent in $\mathfrak{sl}_2$.
Therefore 
\begin{eqnarray*}
\frac{a}{b}, \qquad \frac{1-a}{1-c}, \qquad \frac{d}{c}, \qquad \frac{1-d}{1-b}
\end{eqnarray*}
coincide and this yields $a=b=c=d$ after a brief calculation.
The result follows.
\hfill $\Box $ \\ 

\noindent We finish this section with a comment.
Up to isomorphism
there exists a unique $\boxtimes$-module  with dimension
1, and every element of $\boxtimes$ is 0 on this
$\boxtimes$-module. 
We call this
$\boxtimes$-module
{\it trivial}.

\begin{lemma}
\label{lem:3ind}
Let $V$ denote a nontrivial finite-dimensional irreducible
$\boxtimes$-module. Then for mutually distinct
$i,j,k \in \I$ the actions of $x_{ij}$, $x_{jk}$, $x_{ki}$
on $V$ are linearly independent.
\end{lemma}
\noindent {\it Proof:} 
Let $K=\lbrace \xi \in \boxtimes\,\vert \,\xi.v=0\; \forall v \in V\rbrace$
be the kernel of the $\boxtimes$-action on $V$
and note that $K$ is an ideal of $\boxtimes$. By assumption
$V$ is nontrivial so $K\not=\boxtimes$. Applying Lemma
\ref{lem:ideal} to $K$ we obtain the result.
\hfill $\Box $ \\

\section{Evaluation modules for $\boxtimes$}

\noindent In this section we consider a class of
finite-dimensional irreducible $\boxtimes$-modules
called {\it evaluation modules}.
We start with a lemma.

\begin{lemma}
\label{lem:ev}
With reference to 
Definition 
\ref{def:3point}
the following (i), (ii) hold for
$a \in \K\backslash \lbrace 0,1\rbrace$.
\begin{enumerate}
\item
There exists an $\K$-algebra homomorphism $e_a: \mathcal{A}\to \K$
that sends $t\mapsto a$.
\item There exists a Lie algebra homomorphism ${ev}_a: 
\mathfrak{sl}_2 \otimes \mathcal{A} \to \mathfrak{sl}_2$
that sends $u\otimes b \to u e_a(b)$ for all $u \in \mathfrak{sl}_2$
and $b \in \mathcal{A}$. Moreover 
${ev}_a$ is surjective.
\end{enumerate}
\end{lemma}
\noindent {\it Proof:}  Routine.
\hfill $\Box $ \\ 

\begin{definition}
\label{def:EV}
\rm
For
$a \in \K\backslash \lbrace 0,1\rbrace$
we define
$EV_a : \boxtimes \to \mathfrak{sl}_2$
to be the composition $ev_a \circ \psi$,
where $\psi$ and $ev_a$ are from
 Proposition 
\ref{prop:eta} and Lemma
\ref{lem:ev} respectively.
We note that $EV_a$ is a surjective homomorphism of Lie
algebras.
\end{definition}

\noindent The map $EV_a$ from Definition \ref{def:EV}
is characterized as follows.

\begin{lemma}
\label{lem:evchar}
For a Lie algebra homomorphism $\varepsilon :\boxtimes \to \mathfrak{sl}_2$
the following (i), (ii) are equivalent.
\begin{enumerate}
\item $\varepsilon$ sends
\begin{eqnarray}
\label{eq:esend}
x_{12} \mapsto x,\qquad \qquad 
x_{23} \mapsto y, \qquad \qquad 
x_{31} \mapsto z,
\end{eqnarray}
where $x,y,z$ is the equitable basis for
$\mathfrak{sl}_2$.
\item There exists
$a \in \K\backslash \lbrace 0,1\rbrace$
such that
$\varepsilon=EV_a$.
\end{enumerate}
Suppose (i), (ii) hold. Then $\varepsilon$ sends
\begin{eqnarray}
x_{03} &\mapsto& ay+(a-1)z,
\label{eq:xs}
\\
x_{01} &\mapsto& (1-a^{-1})z-a^{-1}x,
\label{eq:ys}
\\
x_{02} &\mapsto& (1-a)^{-1}x+a(1-a)^{-1}y.
\label{eq:zs}
\end{eqnarray}
\end{lemma}
\noindent {\it Proof:} 
(i) $\Rightarrow$ (ii) 
The map $\varepsilon$ is nonzero so Proposition
\ref{prop:a} applies; let $a \in \K$  be from
that proposition and note that $a\not=0$, $a\not=1$.
We show that $\varepsilon=EV_a$.
To this end we first show that $\varepsilon$ satisfies 
(\ref{eq:xs})--(\ref{eq:zs}). 
Line
 (\ref{eq:xs}) follows from
Definition \ref{def:tet}(i),
line (\ref{eq:esend}),
and the bottom-right display in
 Proposition
\ref{prop:a}.
Line
 (\ref{eq:ys}) follows from
Definition \ref{def:tet}(i),
line (\ref{eq:esend}),
and the top-right display in
 Proposition
\ref{prop:a}.
Line
 (\ref{eq:zs}) follows from
Definition \ref{def:tet}(i),
line (\ref{eq:esend}),
and the bottom-left display in
 Proposition
\ref{prop:a}.
We have now shown that $\varepsilon$ satisfies 
(\ref{eq:xs})--(\ref{eq:zs}). 
Comparing the data in Proposition
\ref{prop:eta} with 
(\ref{eq:esend})--(\ref{eq:zs}) we find
 $\varepsilon(x_{ij})=EV_a(x_{ij})$ for all distinct $i,j \in \I$.
It follows that
$\varepsilon=EV_a$.
\\
\noindent 
(ii) $\Rightarrow$ (i) Immediate from 
Proposition
\ref{prop:eta}
and Definition
\ref{def:EV}.
\\
\noindent Suppose (i), (ii) hold. We saw in the proof of
(i) $\Rightarrow$ (ii) that 
$\varepsilon $ satisfies (\ref{eq:xs})--(\ref{eq:zs}).
\hfill $\Box $ \\

\begin{definition}
\label{def:ev}
\rm
For a finite-dimensional $\mathfrak{sl}_2$-module $V$
 and for
$a \in \K\backslash \lbrace 0,1\rbrace$
 we pull back the 
$\mathfrak{sl}_2$-module action via 
$EV_a$ to obtain a $\boxtimes$-module structure on $V$. We denote
this $\boxtimes$-module by $V(a)$.
\end{definition}

\noindent We now recall the finite-dimensional irreducible
$\mathfrak{sl}_2$-modules.

\begin{lemma}
\label{lem:sl2mod}
{\rm \cite[p. 31]{humphreys}}
There exists a family
\begin{equation}
\label{sl2modlist}
V_d  \qquad \qquad d = 0,1,2,\ldots
\end{equation}
of finite-dimensional irreducible 
$\mathfrak{sl}_2$-modules
with the following
properties. The module $V_d$ has a basis $\lbrace v_n\rbrace_{n=0}^d$
satisfying $h.v_n=(d-2n)v_n$
for $0 \leq n \leq d$,
$f.v_n=(n+1)v_{n+1}$ for $0 \leq n \leq d-1$, $f.v_d=0$,
$e.v_n=(d-n+1)v_{n-1}$ for $1 \leq n \leq d$, $e.v_0=0$.
Every  finite-dimensional irreducible
$\mathfrak{sl}_2$-module
is isomorphic to exactly one of the modules
in line (\ref{sl2modlist}).
\end{lemma}

\begin{lemma}
\label{lem:xyzact}
With respect to the basis $\lbrace v_n \rbrace_{n=0}^d$
for $V_d$ given in 
Lemma
\ref{lem:sl2mod}, the elements $x,y,z$ 
act as follows.
We have $(x+(d-2n)I).v_n= 2(d-n+1)v_{n-1}$
for $1 \leq n \leq d$, 
$(x+dI).v_0=0$,
$(y+(d-2n)I).v_n= -2(n+1)v_{n+1}$
for $0 \leq n \leq d-1$, 
$(y-dI).v_d=0$,
$z.v_n=(d-2n)v_n$ for
$0 \leq n \leq d$.
\end{lemma}
\noindent {\it Proof:} 
Use Lemma
\ref{lem:sl2mod} and
(\ref{eq:xyz}).
\hfill $\Box $ \\

\begin{definition}
\label{def:evalm}
\rm 
With reference to 
Definition \ref{def:ev} and
Lemma \ref{lem:sl2mod},
by an {\it evaluation module} for $\boxtimes$ we mean
a $\boxtimes$-module $V_d(a)$ where $d$ is a positive
integer and
$a \in \K\backslash \lbrace 0,1\rbrace$.
We note that the $\boxtimes$-module $V_d(a)$ is nontrivial and 
irreducible, with diameter $d$ and shape $(1,1,\ldots,1)$. 
We call $a$ the {\it evaluation parameter} of $V_d(a)$.
\end{definition}

\begin{example}
\rm
For $a \in \K\backslash \lbrace 0,1\rbrace$ 
the
$\boxtimes$-module
$V_1(a)$ is described as follows.
Let $v_0, v_1$ denote the basis for the
$\mathfrak{sl}_2$-module $V_1$ given in
Lemma
\ref{lem:sl2mod}.
With respect to this basis and for the $\boxtimes$-module
$V_1(a)$ the generators $x_{ij}$ are represented by the
following matrices.
\begin{eqnarray*}
&&
x_{12}:\;
\left(
\begin{array}{ c c }
-1 & 2   \\
0 & 1  \\
\end{array}
\right), 
\qquad \qquad
\quad \;\;x_{03}:\;
\left(
\begin{array}{ c c }
-1 & 0  \\
-2a & 1   \\
\end{array}
\right), 
\\
&&
x_{23}:\;
\left(
\begin{array}{ c c }
-1 & 0  \\
-2 & 1   \\
\end{array}
\right),
\qquad \qquad
\quad \;\;
x_{01}:\;
\left(
\begin{array}{ c c }
1 & -2a^{-1}   \\
0 & -1  \\
\end{array}
\right), 
\\
&&
x_{31}:\;
\left(
\begin{array}{ c c }
1 & 0   \\
0 & -1  \\
\end{array}
\right), 
\qquad \qquad
\quad \;\;
x_{02}:\;
\left(
\begin{array}{ c c }
\frac{a+1}{a-1} & \frac{2}{1-a}  \\
\frac{2a}{a-1} & \frac{1+a}{1-a} \\
\end{array}
\right) .
\end{eqnarray*}
\end{example}
\noindent {\it Proof:} 
Combine
Lemma
\ref{lem:evchar} and Lemma
\ref{lem:xyzact}.
\hfill $\Box $ \\

\section{The evaluation modules for $\boxtimes$; three characterizations}

\noindent In this section we give three characterizations of
the evaluation modules for $\boxtimes$.
Here is the first characterization.

\begin{proposition}
\label{prop:evchar}
Let $V$ denote a nontrivial finite-dimensional irreducible
$\boxtimes$-module. Then for 
 $a \in \K\backslash \lbrace 0,1\rbrace$ 
the following (i), (ii) are equivalent.
\begin{enumerate}
\item
$V$ is isomorphic to an evaluation module with evaluation parameter $a$.
\item
Each of the following vanishes on $V$:
\begin{eqnarray}
&&
\label{eq:ad1}
a x_{01}+(1-a)x_{02}-x_{03},
\qquad \qquad a x_{10}+(1-a)x_{13}-x_{12},
\\
&&
\label{eq:ad2}
a x_{23}+(1-a)x_{20}-x_{21},
\qquad \qquad 
a x_{32}+(1-a)x_{31}-x_{30}.
\end{eqnarray}
\end{enumerate}
\end{proposition}
\noindent {\it Proof:} 
\noindent (i) $\Rightarrow$ (ii) Immediate from
Lemma
\ref{lem:evchar}.
\\
\noindent 
\noindent (ii) $\Rightarrow$ (i)
By Lemma
\ref{lem:sl2inj}
the generators 
$x_{12}, x_{23}, x_{31}$ form a basis for 
a subalgebra $L$ of $\boxtimes$ that is isomorphic
to $\mathfrak{sl}_2$.
Let $K$ denote the ideal of $\boxtimes$ generated by
the four elements in (\ref{eq:ad1}), (\ref{eq:ad2}).
We claim that $K$ is equal to the 
 kernel of the
$\boxtimes$-action on $V$.
To prove the claim, let $K'$ denote the kernel in
question and observe that $K\subseteq K'$.
Note that $K'$ is an ideal of $\boxtimes$; also
$K'\not=\boxtimes$ since $V$ is nontrivial.
Applying Lemma
\ref{lem:ideal} to $K'$ we obtain 
$K'\cap L=0$. Note  that $K+L$ is a subalgebra
of $\boxtimes$ that contains the generators
$x_{ij}$ $(i,j\in \I,\, i\not=j)$. Therefore
$K+L=\boxtimes$. By these comments
$K=K'$ and
the claim is proved.
We now show that
 the $\boxtimes$-modules $V$ and $V_d(a)$ are isomorphic,
where $d+1$ is the dimension of $V$.
For  $V$ and $V_d(a)$ we restrict the
$\boxtimes$-action to $L$.
We mentioned earlier that the $\boxtimes$-module $V_d(a)$ is
irreducible. Now the $L$-module $V_d(a)$ is irreducible
since $K+L=\boxtimes$ and $K$ vanishes on
$V_d(a)$.
Similarly the $L$-module $V$ is irreducible.
Since $L$ is isomorphic to $\mathfrak{sl}_2$
and since $V$, $V_d(a)$ have the same dimension
the $L$-modules $V$, $V_d(a)$ are isomorphic.
Let $\mu:V\to V_d(a)$ denote an isomorphism of
$L$-modules.
Since $K+L=\boxtimes$ and $K$ vanishes on each 
of $V$, $V_d(a)$ the map $\mu$ 
extends to a $\boxtimes$-module isomomorphism
from $V$ to $V_d(a)$. Therefore the $\boxtimes$-modules
$V$ and $V_d(a)$ are isomorphic as desired.
\hfill $\Box $ \\ 

\noindent We have a comment.

\begin{lemma}
\label{lem:eviso}
Two evaluation modules for $\boxtimes$ are isomorphic
if and only if they have the same diameter and the 
same evaluation parameter.
\end{lemma}
\noindent {\it Proof:} 
Two isomorphic evaluation modules have the
same diameter, since they have the same dimension and the dimension
is one more than the diameter.
Suppose we are given isomorphic evaluation modules
$V_d(a)$ and $V_d(b)$. We show $a=b$.
By 
(\ref{eq:ad1}) the element
$ax_{01}+(1-a)x_{02}-x_{03}$ vanishes on $V_d(a)$.
By 
(\ref{eq:ad1}) and since
$V_d(a), V_d(b)$ are isomorphic the element 
$bx_{01}+(1-b)x_{02}-x_{03}$ vanishes on $V_d(a)$.
Comparing these elements we find that
$a=b$ or  
$x_{01}-x_{02}$ vanishes on $V_d(a)$. The second possibility
 contradicts Definition \ref{def:tet}(i)
and Lemma
\ref{lem:3ind},
so $a=b$ and the result follows.
\hfill $\Box $ \\

\noindent Here is the second characterization of the 
evaluation modules for $\boxtimes$.

\begin{proposition}
\label{prop:evchar2}
Let $V$ denote a nontrivial finite-dimensional irreducible
$\boxtimes$-module.
Then the following (i), (ii) are equivalent.
\begin{enumerate}
\item
$V$ is isomorphic to an evaluation module for $\boxtimes$.
\item
There exist mutually distinct $i,j,k,\ell \in \I$ such
that 
$x_{ij}, 
x_{ik}, 
x_{i\ell}$ are linearly dependent on $V$.
\end{enumerate}
\end{proposition}
\noindent {\it Proof:} 
\noindent (i) $\Rightarrow$ (ii)
Immediate from Proposition
\ref{prop:evchar}.

\noindent (ii) $\Rightarrow$ (i) Without loss
we assume $(i,j,k,\ell)=(0,1,2,3)$.
By assumption there exist $\alpha_1,\alpha_2,\alpha_3 \in \K$,
not all zero, 
such that
\begin{eqnarray}
\label{eq:alphadep}
\alpha_1x_{01}+\alpha_2x_{02}+\alpha_3 x_{03}=0
\end{eqnarray}
on $V$.
Each of 
$\alpha_1,\alpha_2,\alpha_3$ is nonzero by
Definition \ref{def:tet}(i) and
Lemma
\ref{lem:3ind}. Taking the Lie bracket of
$x_{01}$ with 
(\ref{eq:alphadep}) and simplifying the result using
Definition \ref{def:tet}(i),(ii) we find that
\begin{eqnarray}
\label{eq:alphadep2}
(\alpha_2+\alpha_3)x_{01}
-\alpha_2x_{02}
-\alpha_3 x_{03}=0
\end{eqnarray}
on $V$. Adding   
(\ref{eq:alphadep}), 
(\ref{eq:alphadep2}) 
and using
 Lemma
\ref{lem:3ind} we find
\begin{eqnarray}
\label{eq:zsum}
\alpha_1+\alpha_2+\alpha_3=0.
\end{eqnarray}
Taking the Lie bracket of $x_{12}$ with
(\ref{eq:alphadep}) and simplifying the result using
Definition \ref{def:tet}(i),(ii) and
(\ref{eq:zsum})
we find 
\begin{eqnarray}
\label{eq:midj}
\lbrack x_{12}, x_{03}\rbrack &=& 
2\alpha_1\alpha^{-1}_3 x_{01}
-2\alpha_2\alpha^{-1}_3 x_{02}
-2x_{12}
\end{eqnarray}
on $V$. Cyclically permuting $1,2,3$ in the previous argument
we find that both 
\begin{eqnarray}
\label{eq:midk}
\lbrack x_{23}, x_{01} \rbrack &=& 
2\alpha_2\alpha^{-1}_{1}x_{02}
-2\alpha_3\alpha^{-1}_1x_{03}
-2x_{23},
\\
\label{eq:midell}
\lbrack x_{31}, x_{02} \rbrack &=& 
2\alpha_3\alpha^{-1}_2x_{03}
-2\alpha_1\alpha^{-1}_2x_{01}
-2x_{31}
\end{eqnarray}
on $V$.
By the Jacobi identity
\begin{eqnarray}
\lbrack x_{01}, \lbrack x_{12},x_{23}\rbrack \rbrack
+
\lbrack x_{12}, \lbrack x_{23},x_{01}\rbrack \rbrack
+
\lbrack x_{23}, \lbrack x_{01},x_{12}\rbrack \rbrack
=0.
\label{eq:jac}
\end{eqnarray}
Evaluating 
(\ref{eq:jac}) using
Definition \ref{def:tet}(i),(ii) and
(\ref{eq:zsum})--(\ref{eq:midk}) 
we find that
\begin{eqnarray}
\alpha_1 x_{23}+\alpha_2 x_{20} +\alpha_3x_{21} &=& 0
\label{eq:nj}
\end{eqnarray}
on $V$. Cyclically permuting $1,2,3$ in the previous
argument we find that both
\begin{eqnarray}
\alpha_2 x_{31}+\alpha_3 x_{30} +\alpha_1 x_{32} &=& 0,
\label{eq:nk}
\\
\alpha_3 x_{12}+\alpha_1 x_{10} +\alpha_2 x_{13} &=& 0
\label{eq:nell}
\end{eqnarray}
on $V$. By  
(\ref{eq:alphadep}),
(\ref{eq:zsum}) and
(\ref{eq:nj})--(\ref{eq:nell}) 
the four expressions
(\ref{eq:ad1}),
(\ref{eq:ad2}) vanish on 
$V$, where $a=-\alpha_1\alpha^{-1}_{3}$.
Note that $a\not=0$ since $\alpha_1\not=0$, and
$a\not=1$ since 
$\alpha_2 \not=0$.
Now by Proposition
\ref{prop:evchar}, $V$ is isomorphic to an evaluation
module for $\boxtimes$ with evaluation parameter $a$.
The result follows.
\hfill $\Box $ \\ 

\noindent Here is the third characterization of
the evaluation modules for $\boxtimes$.

\begin{proposition}
\label{prop:evchar3}
Let $V$ denote a nontrivial finite-dimensional irreducible
$\boxtimes$-module.
Then the following (i), (ii) are equivalent.
\begin{enumerate}
\item
$V$ is isomorphic to an evaluation module for $\boxtimes$.
\item
$V$ has shape $(1,1,\ldots, 1)$.
\end{enumerate}
\end{proposition}
\noindent {\it Proof:} 
(i) $\Rightarrow$ (ii) 
By Definition
\ref{def:evalm} each evaluation module for $\boxtimes$ has
shape $(1,1,\ldots, 1)$.
\\
\noindent
(ii) $\Rightarrow$ (i)
Throughout this proof  the decomposition 
$\lbrack 1,3\rbrack$ of $V$  will be denoted
by
$\lbrace V_n\rbrace_{n=0}^d$.
By definition, for $0 \leq n \leq d$ the space
$V_n$ is an eigenspace for $x_{13}$
with eigenvalue $2n-d$.
By this and since $x_{13}+x_{31}=0$, the space
$V_n$ is an eigenspace for $x_{31}$
with eigenvalue $d-2n$.
Also $V_n$ has dimension 1 by our
 assumption on the shape.
We will need the action of $x_{30}$ on 
$\lbrace V_n\rbrace_{n=0}^d$.
By the table in Section 3,
\begin{eqnarray}
\label{eq:x30act}
(x_{30}-(d-2n)I)V_n\subseteq V_{n+1}
\qquad \qquad 
(0 \leq n \leq d),
\end{eqnarray}
where we recall $V_{d+1}=0$.
By Lemma
\ref{lem:sl2inj} there exists an injective homomorphism
of Lie algebras from $\mathfrak{sl}_2$ to  $\boxtimes$ that
sends $x,y,z$ to
$x_{12}$,
$x_{23}$,
$x_{31}$ respectively.
Using this homomorphism we pull
back the $\boxtimes$-module structure on $V$ to obtain
an $\mathfrak{sl}_2$-module structure on $V$.
Since $z$ and $x_{31}$ agree on $V$ we find
that for $0 \leq n \leq d$, $V_n$ is an eigenspace for $z$ with eigenvalue 
$d-2n$.
These eigenspaces all have dimension 1 so
the  $\mathfrak{sl}_2$-module $V$
is irreducible. Therefore 
 the $\mathfrak{sl}_2$-module $V$ is isomorphic to
 the $\mathfrak{sl}_2$-module $V_d$
from
Lemma
\ref{lem:sl2mod}.
Let $\lbrace v_n\rbrace_{n=0}^d$ denote the
 basis for $V$ given in Lemma
\ref{lem:sl2mod}.   
By Lemma
\ref{lem:xyzact} we find that for
$0 \leq n \leq d$ the vector
$v_n$ is an eigenvector
for $z$ with eigenvalue $d-2n$; therefore
$v_n$ is a basis for $V_n$.
By Lemma
\ref{lem:xyzact}
and since $y,x_{23}$ agree on $V$,
we find
$(x_{23}+(d-2n)I).v_n= -2(n+1)v_{n+1}$
for $0 \leq n \leq d-1$ and
$(x_{23}-dI).v_d=0$.
By 
(\ref{eq:x30act})
there exist scalars 
$\lbrace \alpha_n\rbrace_{n=1}^d$ in $\K$ such that
$(x_{30}-(d-2n)I).v_n=\alpha_{n+1} v_{n+1}$ for $0 \leq n \leq d-1$
and 
$(x_{30}+dI).v_d=0$.
By Definition \ref{def:tet}(ii),
\begin{eqnarray*}
\lbrack x_{23},x_{30}\rbrack = 2 x_{23}+2x_{30}.
\end{eqnarray*}
For $0 \leq n \leq d$ we apply each side of
this equation to $v_n$
and evaluate the result using the above comments.
The resulting equations show 
that $\alpha_n/n$ is independent of $n$ for $1 \leq n \leq d$.
Denoting this common value by $2a$ we find
$(x_{30}-(d-2n)I).v_n=2a(n+1)v_{n+1}$
for $0 \leq n \leq d-1$.
From our comments so far and since $x_{23}+x_{32}=0$, the element
\begin{eqnarray*}
ax_{32}+(1-a)x_{31}-x_{30}
\end{eqnarray*}
vanishes on each of 
 $\lbrace v_n\rbrace_{n=0}^d$ and hence on $V$.
In particular the actions of $x_{30}$, $x_{31}$, $x_{32}$
on $V$ are linearly dependent.
Invoking Proposition \ref{prop:evchar2} we find
$V$ is isomorphic to
an evaluation module for $\boxtimes$.
\hfill $\Box $ \\

\section{Some automorphisms of $\boxtimes$}

\noindent In this section we consider
how the
symmetric group $S_4$ acts on $\boxtimes$ as a group
of automorphisms.   
We investigate what happens when we 
twist an evaluation module for $\boxtimes$ via
an element of $S_4$.

\medskip
\noindent
For the rest of this paper we identify $S_4$ with the 
group of permutations of $\I$.
We use the cycle notation; for example
$(1,2,3)$ denotes the element of $S_4$ that
sends
 $1\mapsto 2
 \mapsto 3
 \mapsto 1$ and
$0\mapsto 0$.
Note that $S_4$ acts on the set of generators for
$\boxtimes$ by 
permuting the indices:
\begin{eqnarray*}
\sigma(x_{ij}) = x_{\sigma(i),\sigma(j)}
\qquad \qquad \sigma \in S_4, \qquad i,j\in \I, \quad i\not=j.
\end{eqnarray*}
This action leaves invariant the defining relations for $\boxtimes$
and therefore induces an
 action of $S_4$ on $\boxtimes$ as a group
of automorphisms.
This $S_4$-action on $\boxtimes$ induces
an action of $S_4$ 
on the set of $\boxtimes$-modules, as we now explain.

\begin{definition}
\label{def:twist}
\rm
Let $V$ denote a $\boxtimes$-module. For $\sigma \in S_4$ 
there exists a $\boxtimes$-module structure on $V$, called
{\it $V$ twisted via $\sigma$}, that behaves as follows:
for all $\xi \in \boxtimes$ and $v \in V$, 
the vector $\xi.v$ computed in
$V$ twisted via $\sigma$ coincides with the vector
$\sigma^{-1}(\xi).v$ computed in
the original $\boxtimes$-module  $V$.
Sometimes we abbreviate ${}^\sigma V$ for $V$ twisted via $\sigma$.
Observe that $S_4$ acts on the set
of $\boxtimes$-modules, with
$\sigma $ sending
$V$ to ${}^\sigma V$ for all $\sigma \in S_4$
and all $\boxtimes$-modules $V$.
\end{definition}

\noindent For the following  three lemmas 
the proofs are routine and left to the reader.

\begin{lemma} 
\label{lem:twist1}
Let $V$ denote a finite-dimensional irreducible
$\boxtimes$-module.
For distinct $i,j\in \I$ and for $\sigma \in S_4$ the following
coincide:
\begin{enumerate}
\item The decomposition $\lbrack i,j\rbrack$ of $V$;
\item The decomposition $\lbrack \sigma(i),\sigma(j)\rbrack$ of $V$ twisted
via $\sigma$.
\end{enumerate}
\end{lemma}

\begin{lemma} 
\label{lem:twist2}
Let $V$ denote a finite-dimensional irreducible
$\boxtimes$-module.
For $i\in \I$ and for $\sigma \in S_4$ the following
coincide:
\begin{enumerate}
\item The flag $\lbrack i\rbrack$ of $V$;
\item The flag $\lbrack \sigma(i)\rbrack$ of $V$ twisted
via $\sigma$.
\end{enumerate}
\end{lemma}

\begin{lemma} 
\label{lem:twist3}
Let $V$ denote a finite-dimensional irreducible
$\boxtimes$-module.
For  $\sigma \in S_4$ the following
coincide:
\begin{enumerate}
\item The shape of  $V$;
\item The shape of  $V$ twisted
via $\sigma$.
\end{enumerate}
\end{lemma}

\begin{corollary}
\label{lem:twist4}
Let $V$ denote an evaluation module for $\boxtimes$, and pick
$\sigma \in S_4$. Then $V$ twisted
via $\sigma$ is isomorphic to
an evaluation module for $\boxtimes$.
\end{corollary}
\noindent {\it Proof:} Combine Proposition
\ref{prop:evchar3}
and 
Lemma \ref{lem:twist3}.
\hfill $\Box $ \\ 

\noindent For each integer  $d\geq 1$ let $\Delta_d$
denote the set of isomorphism classes of 
 evaluation modules for $\boxtimes$ that have
diameter $d$. 
 By
Definition 
\ref{def:evalm}
and
Lemma
\ref{lem:eviso} the map $a\mapsto V_d(a)$ induces
a bijection from
$\F \backslash \lbrace 0,1\rbrace$ to $\Delta_d$.
Via this bijection 
the $S_4$-action on $\Delta_d$ from Corollary
\ref{lem:twist4}
induces an $S_4$-action on
$\F \backslash \lbrace 0,1\rbrace$. This action
is described in the following lemma and subsequent theorem.

\begin{lemma}
\label{lem:s4f}
There exists an $S_4$-action on the set
$\K\backslash \lbrace 0,1\rbrace$
that does the following.
For all $a \in
\K\backslash \lbrace 0,1\rbrace$,
\begin{enumerate}
\item
$(2,0)$ sends $a\mapsto a^{-1}$;
\item
$(0,1)$ sends $a\mapsto a(a-1)^{-1}$;
\item
$(1,3)$ sends $a \mapsto a^{-1}$.
\end{enumerate}
\end{lemma}
\noindent {\it Proof:} 
Consider the following maps on 
$\K\backslash \lbrace 0,1\rbrace$:
\begin{eqnarray*}
\sigma_1: \; a \mapsto a^{-1},
\qquad \qquad 
\sigma_2: \; a \mapsto a(a-1)^{-1},
\qquad \qquad 
\sigma_3: \; a \mapsto a^{-1}.
\end{eqnarray*}
One routinely verifies that
$\sigma^2_1=1$,
$\sigma^2_2=1$,
$\sigma^2_3=1$ and that
$(\sigma_1 \sigma_2)^3=1$,
$(\sigma_2 \sigma_3)^3=1$,
$(\sigma_1 \sigma_3)^2=1$.
The above equations are the defining relations for
 a well-known presentation for
$S_4$ \cite[p.~105]{humphreys2} and the result follows.
\hfill $\Box $ \\ 

\begin{theorem}
\label{thm:evtwist}
For an integer $d\geq 1$, for $\sigma \in S_4$, and
for $a \in \K\backslash \lbrace 0,1\rbrace$ 
 the following 
are isomorphic:
\begin{enumerate}
\item the $\boxtimes$-module $V_d(a)$ twisted via $\sigma$;
\item the $\boxtimes$-module $V_d(\sigma(a))$.
\end{enumerate}
\end{theorem}
\noindent {\it Proof:}   
We abbreviate $W=V_d(a)$.
% ${}^{\sigma}V_d(a)$ for $V_d(a)$ twisted via $\sigma$.
Without loss we assume that $\sigma$ is one of
$(2,0)$, $(0,1)$, $(1,3)$.
In each case we verify using Proposition
\ref{prop:evchar}  that each of the following vanishes on
$W$:
%$V_d(a)$:
\begin{eqnarray*}
&&
\sigma(a) x_{\sigma(0),\sigma(1)}
+(1-\sigma(a))x_{\sigma(0),\sigma(2)}
-x_{\sigma(0),\sigma(3)},
\\
&&\sigma(a) x_{\sigma(1),\sigma(0)}
+(1-\sigma(a))x_{\sigma(1),\sigma(3)}
-x_{\sigma(1),\sigma(2)},
\\
&&\sigma(a) x_{\sigma(2),\sigma(3)}
+(1-\sigma(a))x_{\sigma(2),\sigma(0)}
-x_{\sigma(2),\sigma(1)},
\\
&&\sigma(a) x_{\sigma(3),\sigma(2)}
+(1-\sigma(a))x_{\sigma(3),\sigma(1)}
-x_{\sigma(3),\sigma(0)}.
\end{eqnarray*}
Now by Definition
\ref{def:twist} and $\sigma^2=1$,
each of the 
following vanishes on 
${}^\sigma W$:
%${}^\sigma (V_d(a))$:
\begin{eqnarray*}
&&
\label{eq:ad1t}
\sigma(a) x_{01}+(1-\sigma(a))x_{02}-x_{03},
\qquad \qquad \sigma(a) x_{10}+(1-\sigma(a))x_{13}-x_{12},
\\
&&
\label{eq:ad2t}
\sigma(a) x_{23}+(1-\sigma(a))x_{20}-x_{21},
\qquad \qquad 
\sigma(a) x_{32}+(1-\sigma(a))x_{31}-x_{30}.
\end{eqnarray*}
 Now by Proposition
\ref{prop:evchar}
the $\boxtimes$-module
${}^\sigma W$ 
%${}^\sigma V_d(a)$ 
has evaluation parameter $\sigma(a)$
and is therefore isomorphic to
$V_d(\sigma(a))$.
\hfill $\Box $ \\ 

\section{The group $G$}

\noindent In Lemma
\ref{lem:s4f}
we gave an action of $S_4$ on
$\K\backslash \lbrace 0,1\rbrace$. In this section
we find the
kernel of this action. We will use the following
notation.

\begin{definition}
\label{def:G}
\rm A partition of $\I$ into two disjoint sets, each with
two elements, is said to have shape $(2,2)$.
Let $P$ denote the set of partitions of $\I$ that
have shape $(2,2)$, and note that $P$ has cardinality 3.
The action of $S_4$ on $\I$ induces an action of $S_4$
on $P$ and this yields
a surjective homomorphism of groups $S_4\to S_3$.
Let $G$ denote the kernel of this homomorphism.
We note that $G$ has cardinality 4 and consists of
\begin{eqnarray*}
(01)(23),
\qquad \qquad 
(02)(13),
\qquad \qquad 
(03)(12)
\end{eqnarray*}
together with the identity element. Observe that $G$ is isomorphic
to the Klein 4-group $\Z_2 \times \Z_2$.
\end{definition}

\begin{lemma}
\label{thm:n1}
With reference to Definition \ref{lem:s4f}, for
$\sigma \in S_4$ the following (i), (ii) are equivalent:
\begin{enumerate}
\item
$\sigma \in G$;
\item $\sigma(a)=a$ for all $a \in \K\backslash \lbrace 0,1\rbrace $.
\end{enumerate}
\end{lemma}
\noindent {\it Proof:} 
Let $G'$ denote the set of elements in
$S_4$ that fix each element of 
$\K\backslash \lbrace 0,1\rbrace $. We show $G'=G$.
The set $\lbrace 3, 1/3,-2,-1/2,3/2,2/3\rbrace $
is
an orbit of $S_4$; therefore
the index $\lbrack S_4: G'\rbrack\geq 6$ so
$G'$ has at most 4 elements.
Also $G'$ contains
$(2,0)(1,3)$
by Lemma
\ref{lem:s4f} and $G'$ is normal in $S_4$
so $G\subseteq G'$. It follows that $G'=G$.
\hfill $\Box $ \\

\begin{corollary}
\label{thm:n2}
For an integer $d\geq 1$,
for $\sigma \in G$, and for
 $a \in \K\backslash \lbrace 0,1\rbrace$ 
 the following
are isomorphic:
\begin{enumerate}
\item the $\boxtimes$-module $V_d(a)$ twisted via $\sigma$;
\item the $\boxtimes$-module $V_d(a)$.
\end{enumerate}
\end{corollary}
\noindent {\it Proof:} 
Combine 
Theorem
\ref{thm:evtwist}
and
Lemma \ref{thm:n1}.
\hfill $\Box $ \\

\section{The relatives of an evaluation parameter}

\noindent In 
Lemma 
\ref{lem:s4f}
we considered an $S_4$-action on the set
$\K\backslash\lbrace 0,1\rbrace$. In this section
we describe the orbits of this action.
We start with a definition.

\begin{definition}
\label{def:rel}
\rm
 Pick  
 $a \in \K\backslash \lbrace 0,1\rbrace$ 
and
mutually distinct $i,j,k,\ell \in \I$.
By the {\it $(i,j,k,\ell)$-relative of $a$} we mean
the scalar $\sigma(a)$ where $\sigma \in S_4$ satisfies
$\sigma(i)=2$,
$\sigma(j)=0$,
$\sigma(k)=1$,
$\sigma(\ell)=3$.
\end{definition}

\begin{lemma}
\label{lem:adjrel}
Pick 
 $a \in \K\backslash \lbrace 0,1\rbrace$ 
and 
mutually distinct $i,j,k,\ell \in \I$.
Let $\alpha$ denote the 
 $(i,j,k,\ell)$-relative of $a$.
Then (i)--(iii) hold below.
\begin{enumerate}
\item
$\alpha^{-1}$ is the $(j,i,k,\ell)$-relative of $a$;
\item
$\alpha(\alpha-1)^{-1}$ is the $(i,k,j,\ell)$-relative of $a$;
\item
$\alpha^{-1}$ is the $(i,j,\ell,k)$-relative of $a$.
\end{enumerate}
\end{lemma}
\noindent {\it Proof:} 
By Definition
\ref{def:rel} there exists
$\sigma \in S_4$ 
that sends the sequence
$(a; i,j,k,\ell)$ to
$(\alpha; 2,0,1,3)$.
To obtain (i), note that
$(2,0)\sigma$ sends
$(a;j,i,k,\ell)$ to $(\alpha^{-1}; 2,0,1,3)$.
To obtain (ii), note that
$(0,1)
\sigma$ sends
$(a; i,k,j,\ell)$ to
$(
\alpha(\alpha-1)^{-1};
2,0,1,3)$.
To obtain (iii), note that
$(1,3)\sigma$
sends
$(a;i,j,\ell,k)$ to $(\alpha^{-1};2,0,1,3)$.
\hfill $\Box $ \\ 

%\noindent Here is another way to view the relatives.
\noindent In the following lemma we interpret 
each relative as a cross ratio (see \cite[p.~48]{conway} for
the details).

\begin{lemma}
\label{lem:bar}
 For
 $a \in \K\backslash \lbrace 0,1\rbrace$ 
 and
mutually distinct $i,j,k,\ell \in \I$
the following (i), (ii) coincide:
\begin{enumerate}
\item
the $(i,j,k,\ell)$-relative of
$a$;
\item
the scalar $\frac{\hat{i}-\hat{\ell}}{\hat{i}-\hat{k}}
\,
 \frac{\hat{j}-\hat{k}}{\hat{j}-\hat{\ell}}
$
where we define $\hat{0}=a$,
$\hat{1}=0$,
$\hat{2}=1$,
$\hat{3}=\infty$.
\end{enumerate}

\end{lemma}
\noindent {\it Proof:} 
Let us call the scalar in (ii) the
{\it $(i,j,k,\ell)$-partner} of $a$. Denoting this
partner by $\alpha$ we observe:
\begin{itemize}
\item
$\alpha^{-1}$ is the $(j,i,k,\ell)$-partner of $a$;
\item
$\alpha(\alpha-1)^{-1}$ is the $(i,k,j,\ell)$-partner of $a$;
\item
$\alpha^{-1}$ is the $(i,j,\ell,k)$-partner of $a$.
\end{itemize}
Note that
$a$ is the $(2,0,1,3)$-partner of $a$ as well as the
$(2,0,1,3)$-relative of $a$.
By these comments 
and Lemma
\ref{lem:adjrel} we find that 
the partner function and the relative
function satisfy the same recursion and same initial
condition. Therefore these functions coincide and the result follows.
\hfill $\Box $ \\

\begin{proposition}
\label{lem:table}
 Pick 
 $a \in \K\backslash \lbrace 0,1\rbrace$ 
 and
mutually distinct $i,j,k,\ell \in \I$.
Then the 
 $(i,j,k,\ell)$-relative of $a$ is given in the
 following table.

\bigskip
\centerline{
\begin{tabular}[t]{c|c}
       {\rm $(i,j,k,\ell)$} & {\rm $(i,j,k,\ell)$-relative }
 \\ \hline  \hline
        $(2,0,1,3)  \quad (0,2,3,1) \quad (1,3,2,0) \quad (3,1,0,2)$
	& $a$
	\\
        $(0,2,1,3)  \quad (2,0,3,1) \quad (1,3,0,2) \quad (3,1,2,0)$
	& $a^{-1}$
	\\
        $(1,0,2,3)  \quad (0,1,3,2) \quad (2,3,1,0) \quad (3,2,0,1)$
	& $1-a$
	\\
        $(0,1,2,3)  \quad (1,0,3,2) \quad (2,3,0,1) \quad (3,2,1,0)$
	& $(1-a)^{-1}$
	\\
        $(2,1,0,3)  \quad (1,2,3,0) \quad (0,3,2,1) \quad (3,0,1,2)$
	& $a(a-1)^{-1}$
	\\
        $(1,2,0,3)  \quad (2,1,3,0) \quad (0,3,1,2) \quad (3,0,2,1)$
	& $1-a^{-1}$
	\end{tabular}}

\bigskip
\end{proposition}
\noindent {\it Proof:} 
Routine calculation using
Lemma
\ref{lem:adjrel} or
Lemma \ref{lem:bar}.
\hfill $\Box $ \\

\section{24 bases for an evaluation module}

\noindent Let $V$ denote an evaluation module for
$\boxtimes$. In this section we display 24 bases for 
$V$ that we find attractive. We compute the action
of the generators $x_{ij}$ on each of these bases.
We start with a definition.

\begin{definition}
\label{def:etapre}
\rm
Let $V$ denote an evaluation module for $\boxtimes$
and pick mutually distinct $i,j,k,\ell\in \I$.
A basis $\lbrace u_n\rbrace_{n=0}^d$ of $V$
is called an {\it $\lbrack i,j,k,\ell\rbrack$-basis}
whenever (i), (ii) hold below:
\begin{enumerate}
\item for $0 \leq n \leq d$ the vector
$u_n$ is contained in component
$n$ of the decomposition $\lbrack k,\ell\rbrack$ for $V$;
\item $\sum_{n=0}^d u_n$ is contained in
component 0 of the flag $\lbrack i \rbrack$ for $V$. 
\end{enumerate}
\end{definition}

\noindent Referring to Definition
\ref{def:etapre}, we will show that $V$ has an $\lbrack i,j,k,\ell\rbrack$-basis
after a few comments.

\begin{lemma}
\label{lem:cequiv}
Let $V$ denote an evaluation module for $\boxtimes$
and pick 
 mutually distinct $i,j,k,\ell \in \I$.
Let $\lbrace u_n\rbrace_{n=0}^d$ denote an
$\lbrack i,j,k,\ell \rbrack$-basis for $V$,
and for $0 \leq n \leq d$ let $u'_n$ denote
any vector in $V$. Then the following are equivalent:
\begin{enumerate}
\item the sequence
$\lbrace u'_n\rbrace_{n=0}^d$ is an
$\lbrack i,j,k,\ell \rbrack$-basis for $V$;
\item there exists a nonzero $\beta \in \K$ such 
that $u'_n = \beta u_n$ for $0 \leq n \leq d$.
\end{enumerate}
\end{lemma}
\noindent {\it Proof:} 
(i) $\Rightarrow$ (ii)  Recall that $V$ has shape $(1,1,\ldots,1)$.
\\
(ii) $\Rightarrow$ (i) Immediate from Definition
\ref{def:etapre}.
\hfill $\Box$ \\ 

\begin{lemma}
\label{lem:sl2013}
\rm
Let $V_d(a)$ denote an evaluation 
module for $\boxtimes$.
Then the basis $\lbrace v_n\rbrace_{n=0}^d$
for the $\mathfrak{sl}_2$-module $V_d$ from
Definition
\ref{lem:sl2mod} is
a
$\lbrack 2,0,1,3\rbrack$-basis for $V_d(a)$.
\end{lemma}
\noindent {\it Proof:} 
By 
(\ref{eq:esend})  and since $x_{13}=-x_{31}$
we find $z+x_{13}$ vanishes on $V_d(a)$.
Now by 
 Lemma
\ref{lem:xyzact}, 
for $0 \leq n \leq d$ we have 
$x_{13}.v_n=(2n-d)v_n$ so
  $v_n$ is contained in component $n$
of the decomposition $\lbrack 1,3\rbrack$ for $V_d(a)$.
Define $\eta=\sum_{n=0}^d v_n$.
Using
(\ref{eq:esend}) 
and Lemma
\ref{lem:xyzact} we find
$x_{23}.\eta=-d\eta$; therefore 
$\eta $ is contained in component $0$ of
the flag $\lbrack 2 \rbrack$ for $V_d(a)$.
Now the basis $\lbrace v_n\rbrace_{n=0}^d$ is a
$\lbrack 2,0,1,3\rbrack$-basis by
Definition
\ref{def:etapre}.
\hfill $\Box$ \\

\begin{lemma}
\label{lem:same}
Let $V$ denote an evaluation module for $\boxtimes$
and pick 
 mutually distinct $i,j,k,\ell \in \I$.
Then for $\sigma \in S_4$ the following are the same:
\begin{enumerate}
\item
 an $\lbrack i,j,k,\ell\rbrack$-basis for $V$;
\item
 a $\lbrack \sigma(i),\sigma(j),\sigma(k),\sigma(\ell)\rbrack$-basis 
 for $V$ twisted via $\sigma$.
\end{enumerate}
\end{lemma}
\noindent {\it Proof:} 
Use
Lemma
\ref{lem:twist1} and
Lemma
\ref{lem:twist2}.
\hfill $\Box$ \\

\begin{proposition}
\label{lem:basis} 
Let $V$ denote an evaluation module
for $\boxtimes$ and
pick mutually distinct $i,j,k,\ell\in \I$.
Then there exists an
$\lbrack i,j,k,\ell\rbrack$-basis for $V$.
\end{proposition}
\noindent {\it Proof:} 
Let $\sigma$ denote the element
of $S_4$ that sends the sequence
$(i,j,k,\ell)$ to $(2,0,1,3)$.
The $\boxtimes$-module ${}^\sigma V$ is isomorphic to an evaluation module
by Corollary 
\ref{lem:twist4},
so 
${}^\sigma V$ 
 has a $\lbrack 2,0,1,3 \rbrack$-basis by Lemma
\ref{lem:sl2013}. Now
$V$ has an $\lbrack i,j,k,\ell\rbrack$-basis
in view of Lemma
\ref{lem:same}.
\hfill $\Box$ \\ 

\noindent 
We now find the action of each generator
$x_{rs}$
on the bases given in Proposition
\ref{lem:basis}. 
To describe these actions we use the following notation.
For an integer $d\geq 0$ 
let $\mbox{Mat}_{d+1}(\K)$ denote the $\K$-algebra
of all $d+1$ by $d+1$ matrices that have entries in $\K$.
We index the rows and columns by $0,1,\ldots, d$.
Let $V$ denote a vector space over $\K$ with dimension
$d+1$ and let 
 $\lbrace u_n\rbrace_{n=0}^d$ denote a basis for $V$. 
For a linear transformation $A:V\to V$
there exists a unique
 $X \in
\mbox{Mat}_{d+1}(\K)$ such that
$A.u_j= \sum_{i=0}^d X_{ij}u_i$ for $0 \leq j\leq d$.
We call $X$ the
{\it matrix that represents $A$ with
respect to $\lbrace u_n\rbrace_{n=0}^d$}.

\begin{lemma}
\label{lem:matrep}
Let $V_d(a)$ denote an evaluation module for
$\boxtimes$. For mutually distinct $i,j,k,\ell \in \I$,
for distinct $r,s\in \I$, and for $\sigma \in S_4$
the following are the same:
\begin{enumerate}
\item the matrix that represents $x_{rs}$ with respect to an
$\lbrack i,j,k,\ell\rbrack$-basis for $V_d(a)$;
\item the matrix that represents 
$x_{\sigma(r),\sigma(s)}$ with respect to a
$\lbrack \sigma(i),\sigma(j),\sigma(k),\sigma(\ell)\rbrack$-basis 
for $V_d(\sigma(a))$.
\end{enumerate}
\end{lemma}
\noindent {\it Proof:} 
Routine using
Definition
\ref{def:twist},
Theorem
\ref{thm:evtwist}, and Lemma
\ref{lem:same}.
\hfill $\Box$ \\

%Let $V=V_d(a)$ denote an evaluation module for $\boxtimes$.
%For mutually distinct
%$i,j,k,\ell \in \I$ and for mutually distinct $r,s \in \I$
%we now find the action of $x_{rs}$
%on a 
%$\lbrack i,j,k,\ell\rbrack$-basis for $V$.
%To describe this action we use the following notation.
%Let $\mbox{Mat}_{d+1}(\K)$ denote the $\K$-algebra
%of all $d+1$ by $d+1$ matrices that have entries in $\K$.
%We index the rows and columns by $0,1,\ldots, d$.
%For  $w \in \boxtimes $ and 
%for a basis $\lbrace u_n\rbrace_{n=0}^d$ 
%of $V$ there exists a unique
% $X \in
%\mbox{Mat}_{d+1}(\K)$ such that
%$w.u_j= \sum_{i=0}^d X_{ij}u_i$ for $0 \leq j\leq d$.
%We call $X$ the
%{\it matrix representing $w$ with
%respect to $\lbrace u_n\rbrace_{n=0}^d$}.

\begin{theorem}
\label{thm:xijact}
Let $V_d(a)$ denote an evaluation module for $\boxtimes$.
For mutually distinct
$i,j,k,\ell \in \I$ and distinct $r,s \in \I$
consider the matrix that represents  $x_{rs}$
with respect to an
$\lbrack i,j,k,\ell\rbrack$-basis for $V_d(a)$.
The entries of this matrix are given in the following table.
All entries not displayed are zero.
\bigskip

\centerline{
\begin{tabular}[t]{c|c c c}
       {\rm generator} & {\rm $(n,n-1)$-entry} & {\rm $(n,n)$-entry} &
       {\rm $(n-1,n)$-entry}
 \\ \hline  \hline
      $x_{\ell k}$ & $0$ & $d-2n$ & $0$
      \\
      $x_{k\ell}$ & $0$ & $2n-d$ & $0$
      \\
 \hline
      $x_{ki}$ & $0$ & $2n-d$ & $2d-2n+2$
      \\
      $x_{ik}$ & $0$ & $d-2n$ & $2n-2d-2$
      \\
     \hline 
      $x_{i\ell}$ & $-2n$ & $2n-d$ & $0$
      \\
      $x_{\ell i}$ & $2n$ & $d-2n$ & $0$
      \\
\hline
      $x_{\ell j}$ & $2\alpha n$ & $d-2n$ & $0$
      \\
      $x_{j \ell}$ & $-2\alpha n$ & $2n-d$ & $0$
      \\
    \hline 
     $x_{jk}$ & $0$ & $d-2n$ & $2(n-d-1)\alpha^{-1}$
      \\
      $x_{kj}$ & $0$ & $2n-d$ & $2(d-n+1)\alpha^{-1}$
      \\
\hline
      $x_{ji}$ & $2\alpha n(\alpha -1)^{-1}$
      & $(d-2n)(\alpha+1)(\alpha-1)^{-1}$ 
      & $2(d-n+1)(1-\alpha)^{-1}$
      \\
      $x_{ij}$ & $2\alpha n(1-\alpha)^{-1}$
      & $(d-2n)(\alpha+1)(1-\alpha)^{-1}$ 
      & $2(d-n+1)(\alpha-1)^{-1}$
\end{tabular}}

\bigskip
In the above table the scalar $\alpha$ denotes the
$(i,j,k,\ell)$-relative of $a$ 
from Definition
\ref{def:rel}.
\medskip
\end{theorem}
\noindent {\it Proof:} 
First consider the special case $(i,j,k,\ell)=(2,0,1,3)$.
Without loss we take our
$\lbrack 2,0,1,3\rbrack$-basis
to be the basis 
 $\lbrace v_n\rbrace_{n=0}^d$  from
Lemma
\ref{lem:sl2013}.
We then obtain the
action of $x_{rs}$ on 
 $\lbrace v_n\rbrace_{n=0}^d$ using
Lemma
\ref{lem:evchar} and
Lemma
\ref{lem:xyzact}. This gives the result
for the special case.
To get the result for
general $(i,j,k,\ell)$,
let  
 $\sigma$  denote the element of $S_4$ 
that sends the sequence $(a;i,j,k,\ell)$
to $(\alpha;2,0,1,3)$.
Using
$\sigma $ and  
 Lemma
\ref{lem:matrep} 
we reduce to  
the above special case. 
\hfill $\Box$ \\

\section{A bilinear form for the evaluation modules}

\noindent Let $V$ denote an evaluation module for $\boxtimes$
and consider the 24 bases for $V$ from 
Definition
\ref{def:etapre}. 
In Section 13 we will compute the
transition matrices between certain pairs of bases among
these 24.
Before we get to this, it is convenient
to introduce a bilinear form on $V$.

\medskip
\noindent We recall a few concepts from linear algebra.
Let $V$ denote a vector space over $\K$ with finite positive
dimension.
By a {\it bilinear form} on $V$ we mean a 
map $\langle \,,\,\rangle : V \times V \to \K$ that satisfies
the following four conditions for all $u,v,w \in V$
and for all $\alpha \in \K$:
(i) $\langle u+v,w \rangle =
 \langle u,w \rangle
 +
 \langle v,w \rangle$;
 (ii)
 $
 \langle \alpha u,v \rangle
 =
 \alpha \langle u,v \rangle
 $;
 (iii)
  $\langle u,v+w \rangle =
   \langle u,v \rangle
   +
   \langle u,w \rangle$;
   (iv)
   $
   \langle u, \alpha v \rangle
   =
   \alpha \langle u,v \rangle
   $.
   We observe
   that a scalar multiple of
   a bilinear form on $V$ is a bilinear form on $V$.
   Let
   $
   \langle\, ,\, \rangle
   $
   denote a bilinear form on $V$.
   This form is said to be {\it symmetric}
   (resp. {\it antisymmetric}) 
  whenever
   $
   \langle u,v \rangle
   =
   \langle v,u \rangle
   $
   (resp. 
   $
   \langle u,v \rangle
   =
   -\langle v,u \rangle
   $)
   for all $u, v \in V$.
   Let
   $
   \langle\, ,\, \rangle
   $
   denote a bilinear form on $V$.
   Then the following
   are equivalent: (i) there exists a nonzero $u \in V$ such
that $
\langle u,v \rangle
= 0
$
for all $v \in V$;
 (ii) there exists a nonzero $v \in V$ such
 that $
 \langle u,v \rangle
 = 0
 $
 for all $u \in V$.
 The form
 $\langle \,,\,\rangle $
 is
 said to be {\it degenerate }
 whenever (i), (ii) hold and {\it nondegenerate}
 otherwise.
Assume 
 $\langle \,,\,\rangle $ is nondegenerate.
Let $\lbrace U_n\rbrace_{n=0}^d $ 
and 
 $\lbrace V_n\rbrace_{n=0}^d $  denote decompositions
 of $V$ that have the same diameter. These decompositions
 are called {\it dual}
with respect to 
 $\langle \,,\,\rangle $ whenever
 $\langle U_i,V_j\rangle =0$ for $0 \leq i,j\leq d, i \not=j$.
In this case 
$\lbrace U_n\rbrace_{n=0}^d $ 
and 
 $\lbrace V_n\rbrace_{n=0}^d $  have the same shape.

\begin{lemma}
\label{lem:bilsl2}
For each integer $d\geq 0$ there exists a nonzero
bilinear form 
 $\langle \,,\,\rangle $ on the 
$\mathfrak{sl}_2$-module $V_d$ such that
\begin{eqnarray*}
\langle \xi.u,v\rangle = - \langle u, \xi.v\rangle
\qquad \qquad \xi \in \mathfrak{sl}_2,\qquad u,v\in V_d.
\end{eqnarray*}
The form is nondegenerate. The form is unique up to multiplication
by a nonzero scalar in $\K$. The form is symmetric (resp. antisymmetric)
when $d$ is even (resp. odd).
\end{lemma}
\noindent {\it Proof:} 
Concerning existence, let $\lbrace v_n\rbrace_{n=0}^d$ denote
the basis for $V_d$ from Lemma
\ref{lem:sl2mod}, and let 
 $\langle \,,\,\rangle $ denote the bilinear form on $V_d$
 that satisfies
\begin{eqnarray*}
 \langle v_r,v_s\rangle = \delta_{r+s,d}(-1)^r {d \choose r}
\qquad \qquad (0 \leq r,s\leq d).
\end{eqnarray*}
Using the data in Lemma
\ref{lem:sl2mod} we find
 $\langle \xi.v_r,v_s\rangle = 
 -\langle v_r,\xi.v_s\rangle$ for 
 $\xi \in \lbrace e,f,h\rbrace$ and
 $0 \leq r,s\leq d$.
  Since
$\lbrace v_n\rbrace_{n=0}^d$ is a basis
 for $V_d$ and $e,f,h$ is a basis for 
$\mathfrak{sl}_2$  we conclude
 $\langle \xi.u,v\rangle = 
 -\langle u,\xi.v\rangle$ for $\xi \in \mathfrak{sl}_2$
and $u,v\in V_d$.
This shows that the required bilinear form exists.
The remaining assertions are routinely verified.
\hfill $\Box$ \\

\begin{lemma}
\label{lem:bil}
Let $V$ denote an evaluation module for $\boxtimes$.
Then there exists a nonzero
bilinear form 
 $\langle \,,\,\rangle $
on $V$ such that
\begin{eqnarray*}
\langle \xi.u,v\rangle = - \langle u, \xi.v\rangle
\qquad \qquad 
\xi \in \boxtimes, \qquad u,v\in V.
\end{eqnarray*}
The form is nondegenerate. The form is unique up to multiplication
by a nonzero scalar in $\K$. The form is symmetric (resp. antisymmetric)
when the diameter of $V$ is even (resp. odd).
\end{lemma}
\noindent {\it Proof:} 
By definition there exist an integer $d\geq 1$ and
$a \in \K\backslash \lbrace 0,1\rbrace$ 
such that
$V=V_d(a)$. Let 
$\langle \,,\,\rangle$ denote the bilinear form
on 
the $\mathfrak{sl}_2$-module $V_d$ from
Lemma
\ref{lem:bilsl2}. This form meets the
requirements of the present lemma since
each element of $\boxtimes$ acts on $V_d(a)$
as an element of
$\mathfrak{sl}_2$. This shows that the required
bilinear form exists. The remaining assertions
follow from
Lemma
\ref{lem:bilsl2} and since 
$EV_a:\boxtimes\to \mathfrak{sl}_2$ is surjective.
\hfill $\Box$ \\ 

\begin{definition}
\rm Referring to Lemma
\ref{lem:bil} 
we call $\langle\,,\,\rangle$ a {\it standard bilinear
form} for $V$.
\end{definition}

\noindent We have a comment.
\begin{lemma}
Let $V$ denote an evaluation module for $\boxtimes$.
Then for $\sigma \in S_4$ the following (i), (ii) are the
same:
\begin{enumerate}
\item  a
 standard bilinear form
for $V$;
\item
a standard bilinear form for  
$V$ twisted via $\sigma$.
\end{enumerate}
\end{lemma}
\noindent {\it Proof:} 
Combine 
Definition \ref{def:twist}
and Lemma
\ref{lem:bil}.
\hfill $\Box$ \\

\begin{lemma}
\label{lem:dual}
Let $V$ denote an evaluation module for $\boxtimes$
with standard bilinear form 
$\langle\,,\,\rangle$.
Then for distinct $i,j\in \I$ the decompositions
$\lbrack i,j\rbrack$ and 
$\lbrack j,i\rbrack$ of $V$
are dual with respect to 
 $\langle \,,\,\rangle $.
\end{lemma}
\noindent {\it Proof:} 
Let $\lbrace V_n\rbrace_{n=0}^d$ denote the decomposition
$\lbrack i,j\rbrack$ and note that
$\lbrace V_{d-n}\rbrace_{n=0}^d$
is the decomposition
$\lbrack j,i\rbrack$.
For $0 \leq r,s\leq d$ with $r+s\not=d$ we show $V_r,V_s$
are orthogonal with respect to
$\langle \,,\,\rangle$. For $u \in V_r$ and $v\in V_s$
we have
$\langle x_{ij}.u,v\rangle = -
\langle u,x_{ij}.v\rangle$; evaluating this
using $x_{ij}.u=(2r-d)u$ and
 $x_{ij}.v=(2s-d)v$ we obtain $(r+s-d) 
\langle u,v\rangle=0$.
We assume $r+s\not=d$ so 
$\langle u,v\rangle=0$ and therefore
 $V_r,V_s$
are orthogonal with respect to
$\langle \,,\,\rangle$. The result follows.
\hfill $\Box$ \\

\begin{lemma}
\label{lem:flagoc}
Let $V_d(a)$ denote an evaluation module for $\boxtimes$
with standard 
bilinear form  
  $\langle\,,\,\rangle$.
Then for $i \in \I$ and for $0 \leq n \leq d-1$ the following
(i), (ii) 
are orthogonal complements with respect to 
 $\langle \,,\,\rangle $:
\begin{enumerate}
\item component $n$ of the flag $\lbrack i \rbrack$;
\item component $d-n-1$ of the flag $\lbrack i \rbrack$.
\end{enumerate}
\end{lemma}
\noindent {\it Proof:} 
Fix $j \in \I$ $(j\not=i)$ and
let $\lbrace V_n\rbrace_{n=0}^d$ denote the
decomposition $\lbrack i,j\rbrack$ of $V_d(a)$.
By construction, component $n$ (resp. component $d-n-1$)
of the flag $\lbrack i \rbrack $ is 
$V_0+V_1+\cdots +V_n$
(resp. $V_0+V_1+\cdots +V_{d-n-1}$).
These components are orthogonal
by Lemma
\ref{lem:dual}. Moreover the sum of their dimensions
is $d+1$ and this equals the dimension of $V_d(a)$.
Therefore these components are orthogonal complements.
\hfill $\Box$ \\ 

\begin{proposition}
\label{lem:fins}
Let $V$ denote an evaluation module for $\boxtimes$
with standard bilinear form 
$\langle\,,\,\rangle$. Pick mutually distinct
$i,j,k,\ell\in \I$ and consider an 
$\lbrack i,j,k,\ell\rbrack$-basis for $V$
from Definition
\ref{def:etapre}. Denoting this basis by
$\lbrace u_n\rbrace_{n=0}^d$ we have
\begin{eqnarray}
\label{eq:fins}
\langle u_r, u_s\rangle = 
\delta_{r+s,d} (-1)^r
{d \choose r} 
\langle u_0, u_d\rangle 
\end{eqnarray}
for $0 \leq r,s\leq d$.
\end{proposition}
\noindent {\it Proof:} 
If $r+s\not=d$ then
$\langle u_r, u_s\rangle = 0$ by Lemma
\ref{lem:dual}.
Also by Lemma 
\ref{lem:bil} we have 
\begin{eqnarray}
\langle x_{ki}.u_n,u_{d-n+1}\rangle
=-\langle u_n,x_{ki}.u_{d-n+1}\rangle
\label{eq:ss}
\end{eqnarray}
for $1 \leq n\leq d$. The action of
$x_{ki}$ on
$\lbrace u_n\rbrace_{n=0}^d$ is given 
in Theorem
\ref{thm:xijact}.
Evaluating
(\ref{eq:ss}) using this data  we find
\begin{eqnarray*}
(d-n+1)
\langle u_{n-1},u_{d-n+1}\rangle
=
-n
\langle u_n,u_{d-n}\rangle
\qquad \qquad (1 \leq n \leq d).
\end{eqnarray*}
Solving this recursion we find
\begin{eqnarray*}
\langle u_r, u_{d-r}\rangle = 
(-1)^r
{d \choose r} 
\langle u_0, u_d\rangle
\qquad \qquad (0 \leq r \leq d).
\end{eqnarray*}
The result follows.
\hfill $\Box$ \\

%%%%%%%%%%%%%%%%%
\section{A normalization for the 24 bases}

\noindent Let $V$ denote an evaluation module for $\boxtimes$
and consider the 24 bases for $V$ from
Definition 
\ref{def:etapre}.
In the next section we will compute the
transition matrices between certain pairs of bases among
these 24.
In order to do this elegantly
we first normalize our bases.

\begin{notation}
\label{not}
\rm
We fix an integer $d\geq 1$ and a scalar
$a \in \K\backslash \lbrace 0,1\rbrace$.
 We consider the $\boxtimes$-module
$V=V_d(a)$ from Definition
\ref{def:evalm}. For $i \in \I$ we fix a nonzero vector
$\eta_i$ in component $0$ 
of the flag $\lbrack i \rbrack$ for $V$.
Let $\langle\,,\,\rangle$ denote a standard bilinear form on $V$.
\end{notation}

\begin{lemma}
\label{def:eta}
With reference to 
Notation \ref{not}, for 
mutually distinct $i,j,k,\ell \in \I$ there exists a unique
basis $\lbrace u_n \rbrace_{n=0}^d$ of $V$ 
such that
(i), (ii) hold below:
\begin{enumerate}
\item for $0 \leq n \leq d$ the vector
$u_n$ is contained in component
$n$ of the decomposition $\lbrack k,\ell\rbrack$ for $V$;
\item $\eta_i = \sum_{n=0}^d u_n$.
\end{enumerate}
We denote this basis by $\lbrack i,j,k,\ell\rbrack$.
\end{lemma}
\noindent {\it Proof:} 
By Proposition  
\ref{lem:basis}  there exists an
$\lbrack i,j,k,\ell\rbrack$-basis for $V$.
Denote this basis by $\lbrace u'_n \rbrace_{n=0}^d$. 
By Definition 
\ref{def:etapre}, for $0 \leq n \leq d$ the
vector $u'_n$ is in component $n$ of the decomposition
$\lbrack k,\ell\rbrack$. Also
$\sum_{n=0}^d u'_n$ lies in component $0$ of the flag $\lbrack i \rbrack$.
Define $\eta=\sum_{n=0}^d u'_n$ and note that $\eta \not=0$.
Each of $\eta, \eta_i$ span component $0$ of
the flag $\lbrack i \rbrack$ so there exists a nonzero $\beta \in \F$
such that $\eta=\beta \eta_i$. Define $u_n = \beta^{-1}u'_n$ for
$0 \leq n \leq d$. By construction $\lbrace u_n\rbrace_{n=0}^d$
is a basis for $V$ that satisfies (i), (ii) above.
This shows that the required basis exists. This basis is
unique by condition (ii) above and Lemma
\ref{lem:cequiv}.
\hfill $\Box$ \\

\begin{lemma}
\label{lem:fudge}
With reference to
Notation \ref{not},
the following (i), (ii) hold.
\begin{enumerate}
\item
$\langle \eta_i,\eta_i\rangle =0$ for $i \in \I$;
\item
$\langle \eta_i,\eta_j\rangle \not=0$ for distinct $i,j \in \I$.
\end{enumerate}
\end{lemma}
\noindent {\it Proof:} 
\noindent (i) The vector $\eta_i$ is in component $0$
of the flag $\lbrack i \rbrack $ on $V$. This component is
orthogonal to itself by Lemma
\ref{lem:flagoc} and since $d\geq 1$.
\\
\noindent (ii) Let $\lbrace V_n \rbrace_{n=0}^d$ denote the
decomposition $\lbrack i,j\rbrack $ of $V$.
By construction $\eta_i$ spans $V_0$ and $\eta_j$ spans $V_d$.
By Lemma
\ref{lem:flagoc} 
the orthogonal complement of $V_0$ is
 $V_0+\cdots + V_{d-1}$.
Therefore $V_0, V_d$
are not orthogonal and this implies $\langle \eta_i,\eta_j\rangle \not=0$.
\hfill $\Box$ \\ 

\begin{lemma}
\label{lem:unorm}
With reference to Notation \ref{not},
pick mutually distinct $i,j,k,\ell \in \I$
and consider the basis 
$\lbrack i,j,k,\ell\rbrack$
for $V$ 
from Lemma 
\ref{def:eta}.
Denoting this basis by $\lbrace u_n\rbrace_{n=0}^d$ 
we have:
\begin{enumerate}
\item
$u_0 = \eta_k \frac{\langle \eta_i, \eta_{\ell} \rangle}
{\langle \eta_k,\eta_{\ell} \rangle}$;
\item
$u_d = \eta_{\ell} \frac{\langle \eta_k, \eta_i \rangle}
{\langle \eta_k,\eta_{\ell} \rangle}$.
\end{enumerate}
\end{lemma}
\noindent {\it Proof:} 
(i) By construction there exists $\gamma \in \K$ such
that $u_0=\gamma \eta_k$. We show
$\gamma = \langle \eta_i, \eta_{\ell} \rangle
\langle \eta_k,\eta_{\ell} \rangle^{-1}$.
By 
Lemma \ref{lem:dual}
we have
 $\langle u_n,\eta_{\ell}\rangle=0$
for $1 \leq n \leq d$.
By this and Lemma
\ref{def:eta}(ii),
\begin{eqnarray*}
\langle \eta_i, \eta_{\ell}\rangle
&=&
\langle u_0+\cdots+u_d, \eta_{\ell}\rangle
\\
&=&
\langle u_0, \eta_{\ell}\rangle
\\
&=& \gamma \langle \eta_k, \eta_{\ell}\rangle.
\end{eqnarray*}
Therefore 
$\gamma = \langle \eta_i, \eta_{\ell} \rangle
\langle \eta_k,\eta_{\ell} \rangle^{-1}$.
\\
(ii) Similar to the proof of (i) above.
\hfill $\Box$ \\ 

\begin{proposition}
\label{lem:dualb}
With reference to Notation \ref{not},
 pick mutually distinct $i,j,k,\ell \in \I$
and consider the 
basis
$\lbrack i,j,k,\ell\rbrack$ for $V$ from
 Lemma 
\ref{def:eta}.
Denoting this basis by $\lbrace u_n\rbrace_{n=0}^d$
we have
\begin{eqnarray}
\label{eq:uform}
\langle u_r, u_s\rangle = 
\delta_{r+s,d} (-1)^r
{d \choose r} \frac{\langle \eta_k,\eta_i\rangle
\langle \eta_i,\eta_\ell\rangle}{
\langle \eta_k,\eta_{\ell}\rangle}
\end{eqnarray}
for $0 \leq r,s\leq d$.
\end{proposition}
\noindent {\it Proof:} 
In line
(\ref{eq:fins}) evaluate
$\langle u_0,u_d\rangle $ using
Lemma \ref{lem:unorm}.
\hfill $\Box$ \\

\section{Transition matrices between the 24 bases}

\noindent In Lemma \ref{def:eta}
we defined 24 bases for $V_d(a)$.
In this section we compute the transition matrices between
certain pairs of bases among these 24.
 First we recall a few terms. 

\medskip
\noindent Let $V$ denote a vector space over $\K$ with finite
positive dimension. 
Suppose we are given two bases for $V$, written
$\lbrace u_n\rbrace_{n=0}^d$ and
$\lbrace v_n\rbrace_{n=0}^d$.
By the {\it transition matrix}
from
$\lbrace u_n\rbrace_{n=0}^d$ to 
$\lbrace v_n\rbrace_{n=0}^d$
 we mean the matrix
$T \in \hbox{Mat}_{d+1}(\K)$ such that
\begin{equation}
v_j = \sum_{i=0}^d T_{ij}u_i      \qquad \qquad (0 \leq j\leq d).
\label{eq:transdefpre}
\end{equation}
We recall a few properties of transition matrices.
Let $T$ denote the transition matrix from
$\lbrace u_n\rbrace_{n=0}^d$ to 
$\lbrace v_n\rbrace_{n=0}^d$.
Then $T^{-1}$ exists, and equals
the transition matrix from
$\lbrace v_n\rbrace_{n=0}^d$ to 
$\lbrace u_n\rbrace_{n=0}^d$.
Let
$\lbrace w_n\rbrace_{n=0}^d$ 
denote a basis for $V$, and let
$S$ denote the transition matrix from
$\lbrace v_n\rbrace_{n=0}^d$ to 
$\lbrace w_n\rbrace_{n=0}^d$.
Then $TS$ is the transition
matrix from
$\lbrace u_n\rbrace_{n=0}^d$ to 
$\lbrace w_n\rbrace_{n=0}^d$.
For a linear transformation $A:V\to V$ let
$X$ (resp. $Y$) denote the matrix
representing $A$ with respect to
$\lbrace u_n\rbrace_{n=0}^d$
(resp. $\lbrace v_n\rbrace_{n=0}^d$).
Then $XT = TY$.

\medskip
\noindent
The following matrix will play a role in our discussion.
For an integer $d\geq 0$
let $Z=Z(d)$ denote the matrix in
$\mbox{Mat}_{d+1}(\K)$ with  entries
\begin{equation}
Z_{ij} = \cases{1, &if $\;i+j=d$;\cr
0, &if $\;i+j \not=d$\cr}
\qquad \qquad
(0 \leq i,j\leq d).
\label{eq:zmat}
\end{equation}
We observe $Z^2=I$.

\begin{theorem}
\label{thm:trans}
With reference to Notation \ref{not} and Lemma 
\ref{def:eta}, pick mutually distinct $i,j,k,\ell \in \I$
and consider the transition matrices from the basis
$\lbrack i,j,k,\ell \rbrack$ to the bases
\begin{eqnarray*}
\lbrack j,i,k,\ell \rbrack,
\qquad \qquad 
\lbrack i,k,j,\ell\rbrack,
\qquad \qquad 
\lbrack i,j,\ell,k\rbrack.
\end{eqnarray*}
\begin{enumerate}
\item The first transition matrix is diagonal with
$(r,r)$-entry
\begin{eqnarray*}
\frac{\langle \eta_j,\eta_{\ell} \rangle}
{\langle \eta_i,\eta_{\ell} \rangle} \,\alpha^r
\end{eqnarray*}
for $0 \leq r \leq d$, where $\alpha$ is the $(i,j,k,\ell)$-relative
of $a$  from Definition
\ref{def:rel}.
\item The second transition matrix is lower triangular with $(r,s)$-entry
\begin{eqnarray*}
{r \choose s}\alpha^{r-s}(1-\alpha)^s
\end{eqnarray*}
for $0 \leq s\leq r \leq d$, 
where $\alpha$ is the $(i,j,k,\ell)$-relative
of $a$  from Definition
\ref{def:rel}.
\item The third transition matrix is the matrix $Z$ from
(\ref{eq:zmat}).
\end{enumerate}
\end{theorem}
\noindent {\it Proof:} 
\noindent (i)
Let $T$ denote the transition matrix
from
$\lbrack i,j,k,\ell \rbrack$ to
$\lbrack j,i,k,\ell \rbrack$.
Then 
$T$ is diagonal by Lemma
\ref{def:eta}(i).
Let $X$ (resp. $Y$) denote the matrix that
represents $x_{j\ell}$ with respect to
$\lbrack i,j,k,\ell \rbrack$ (resp.  
$\lbrack j,i,k,\ell \rbrack$).
Then $XT=TY$.
By Theorem \ref{thm:xijact},
$X$ is lower bidiagonal with entries
$X_{nn}= 2n-d$ for $0 \leq n \leq d$
and $X_{n,n-1}=-2\alpha n$ for
$1 \leq n\leq d$.
Similarly 
$Y$ is lower bidiagonal with entries
$Y_{nn}= 2n-d$ for $0 \leq n \leq d$
and $Y_{n,n-1}=-2n$ for
$1 \leq n\leq d$.
Evaluating $XT=TY$ using these comments we find
$T_{nn}=\alpha T_{n-1,n-1}$ for $1\leq n \leq d$.
Therefore 
\begin{eqnarray}
T_{rr}=T_{00}\alpha^r  \qquad \qquad (0 \leq r \leq d).
\label{eq:tdiag}
\end{eqnarray}
By Lemma
\ref{lem:unorm}(i) the $0$th component of the 
basis 
$\lbrack i,j,k,\ell \rbrack$ 
is 
$\eta_k \langle \eta_i, \eta_{\ell} \rangle
\langle \eta_k,\eta_{\ell} \rangle^{-1}$.
Also by Lemma
\ref{lem:unorm}(i), the $0$th component of the 
basis 
$\lbrack j,i,k,\ell \rbrack$ 
is 
$\eta_k \langle \eta_j, \eta_{\ell} \rangle
\langle \eta_k,\eta_{\ell} \rangle^{-1}$.
By these comments
the $0$th component of the 
basis 
$\lbrack j,i,k,\ell \rbrack$ 
is  
$\langle \eta_j, \eta_{\ell} \rangle
\langle \eta_i,\eta_{\ell} \rangle^{-1}$ times
the $0$th component of the 
basis 
$\lbrack i,j,k,\ell \rbrack$.
Therefore 
$T_{00} = 
\langle \eta_j, \eta_{\ell} \rangle
\langle \eta_i,\eta_{\ell} \rangle^{-1}$. 
Combining this with
(\ref{eq:tdiag}) we get the result.
\\
\noindent (ii) Let $T$ denote the transition matrix
from 
$\lbrack i,j,k,\ell \rbrack$ to
$\lbrack i,k,j,\ell \rbrack$.
For each of these two bases the sum of the vectors is $\eta_i$, so
$T$ has constant row sum 1.
By construction 
$T$ is lower triangular.
Let $X$ (resp. $Y$) denote the matrix
that represents $x_{\ell j}$ with respect to
$\lbrack i,j,k,\ell \rbrack$ 
(resp. 
$\lbrack i,k,j,\ell \rbrack$).
Then $XT=TY$.
By Theorem
\ref{thm:xijact}, $X$ is lower bidiagonal with
entries
$X_{nn}=d-2n$ for $0 \leq n \leq d$ and
$X_{n,n-1}=2\alpha n$ for $1\leq n\leq d$.
Also $Y$ is diagonal with entries
$Y_{nn}=d-2n$ for $0 \leq n \leq d$.
Evaluating $XT=TY$ using these comments we find
$(r-s)T_{rs}=\alpha r T_{r-1,s}$ for $0 \leq s < r \leq d$.
By this recursion,
\begin{eqnarray}
T_{rs} = 
{r \choose s}\alpha^{r-s}T_{ss}
\qquad \qquad (0 \leq s \leq r \leq d). 
\label{eq:tdata}
\end{eqnarray}
By 
(\ref{eq:tdata}) 
and since $T$ has constant row sum 1 we routinely
obtain $T_{ss}=(1-\alpha)^s$
for $0 \leq s \leq d$ by induction on $s$.
The result follows.
\\
\noindent (iii) Immediate from
Lemma \ref{def:eta}.
\hfill $\Box$ \\

\noindent With reference to
 Notation \ref{not}, 
we now consider how the scalars
\begin{eqnarray*}
\langle \eta_i,\eta_j\rangle \qquad \qquad i,j\in \I, \quad i\not=j
\end{eqnarray*}
are related.

\begin{lemma}
\label{lem:f}
With reference to Notation \ref{not},
for mutually distinct $i,j,k,\ell \in \I$ we have
\begin{eqnarray*}
\frac{\langle \eta_i, \eta_\ell \rangle}
{\langle \eta_i, \eta_k \rangle}
\,\frac{\langle \eta_j, \eta_k \rangle}
{\langle \eta_j, \eta_\ell \rangle}
&=& \alpha^d,
\end{eqnarray*}
where $\alpha$ is the $(i,j,k,\ell)$-relative
of $a$ from Definition
\ref{def:rel}.
\end{lemma}
\noindent {\it Proof:} 
Let $T$ (resp. $T'$) 
denote the transition matrix from
$\lbrack i,j,k,\ell\rbrack$
to
$\lbrack j,i,k,\ell\rbrack$
(resp. 
$\lbrack j,i,\ell,k\rbrack$
to
$\lbrack i,j,\ell,k\rbrack$).
By Theorem
\ref{thm:trans}(iii), the matrix
$Z$ is the transition matrix
from
$\lbrack j,i,k,\ell\rbrack$
to
$\lbrack j,i,\ell,k\rbrack$,
and from
$\lbrack i,j,\ell,k\rbrack$
to
$\lbrack i,j,k,\ell\rbrack$.
Therefore $TZT'Z=I$.
By Theorem
\ref{thm:trans}(i),
$T$ is diagonal with entries
$
T_{nn} = 
\langle \eta_j,\eta_{\ell} \rangle
\langle \eta_i,\eta_{\ell} \rangle^{-1} \alpha^n
$
for $0 \leq n \leq d$.
By Lemma
\ref{lem:adjrel}(i),(iii) the
$(j,i,\ell,k)$-relative of $a$
coincides with the 
$(i,j,k,\ell)$-relative of $a$ and is therefore
equal to 
 $\alpha$.
By this 
and
 Theorem 
\ref{thm:trans}(i),
$T'$ is diagonal with entries
$
T'_{nn} = 
\langle \eta_i,\eta_{k} \rangle
\langle \eta_j,\eta_{k} \rangle^{-1} \alpha^n
$
for $0 \leq n \leq d$.
Evaluating 
 $TZT'Z=I$ using these comments we obtain
 the result.
\hfill $\Box$ \\ 

\begin{corollary}
\label{cor:3eq}
With reference to Notation \ref{not} 
we have
\begin{eqnarray*}
\frac{\langle \eta_0, \eta_1 \rangle}
{\langle \eta_0, \eta_3 \rangle}
\,\frac{\langle \eta_2, \eta_3 \rangle}
{\langle \eta_2, \eta_1 \rangle}
&=& a^d,
\\
\frac{\langle \eta_0, \eta_2 \rangle}
{\langle \eta_0, \eta_1 \rangle}
\,\frac{\langle \eta_3, \eta_1 \rangle}
{\langle \eta_3, \eta_2 \rangle}
&=& (1-a^{-1})^d,
\\
\frac{\langle \eta_0, \eta_3 \rangle}
{\langle \eta_0, \eta_2 \rangle}
\,\frac{\langle \eta_1, \eta_2 \rangle}
{\langle \eta_1, \eta_3 \rangle}
&=& (1-a)^{-d}.
\end{eqnarray*}
\end{corollary}
\noindent {\it Proof:} 
Use
Lemma \ref{lem:f} and
Proposition
\ref{lem:table}.
\hfill $\Box$ \\

\begin{note}
\label{note:free}
\rm
By Corollary \ref{cor:3eq} and the symmetry/antisymmetry of
$\langle \,,\,\rangle$ the
scalars
\begin{eqnarray*}
\langle \eta_i, \eta_j\rangle \qquad \qquad i,j\in \I, \quad i\not=j
\end{eqnarray*}
are determined by the sequence
\begin{eqnarray}
\langle \eta_0, \eta_1\rangle,
\qquad
\langle \eta_0, \eta_2\rangle,
\qquad
\langle \eta_0, \eta_3\rangle,
\qquad
\langle \eta_1, \eta_2\rangle.
\qquad
\label{eq:free}
\end{eqnarray}
The scalars
(\ref{eq:free})
are ``free'' in the following sense.
Given a sequence $\Psi$ of four nonzero scalars in $\K$,
there exist vectors 
$\eta_i$ $(i \in \I)$ as in Notation \ref{not}
such that  
the sequence (\ref{eq:free}) is equal to $\Psi$.
\end{note}

\section{Realizing the evaluation modules for $\boxtimes$
using polynomials in two variables}

\noindent In this section we give a concrete realization of
the evaluation modules for $\boxtimes$ using polynomials
in two variables.

\begin{notation}
\label{not:a}
\rm
Let $z_0, z_1$ denote commuting indeterminates.
Let $\K\lbrack z_0,z_1\rbrack$ denote the $\K$-algebra
of all polynomials in $z_0,z_1$ that have coefficients
in $\K$. We abbreviate ${\mathcal P}=
\K\lbrack z_0,z_1\rbrack$ and often view this as a vector
space over $\K$.
For an integer $d\geq 0$ let ${\mathcal P}_d$
denote the subspace of 
${\mathcal P}$ consisting
of the homogeneous polynomials in $z_0,z_1$ that have
total degree $d$. Thus
$\lbrace z^{d-n}_0z^n_1\rbrace_{n=0}^d$ is a basis for
${\mathcal P}_d$.
Note that
\begin{eqnarray}
\label{eq:ads}
{\mathcal P} = \sum_{n=0}^{\infty} {\mathcal P}_d
\qquad \qquad (\mbox{direct sum})
\end{eqnarray}
and that
${\mathcal P}_r {\mathcal P}_s = {\mathcal P}_{r+s}$
for $r,s\geq 0$. 
We fix mutually distinct 
$\beta_i \in \K$ $(i \in \I)$.
Then there exist unique $z_2,z_3 \in {\mathcal P}$ such that
\begin{eqnarray}
\sum_{i\in \I} z_i = 0,
\qquad \qquad 
\sum_{i\in \I} \beta_i z_i = 0.
\label{eq:linzz}
\end{eqnarray}
\end{notation}

\noindent For the following three lemmas the proofs
are routine and left to the reader.

\begin{lemma}
\label{lem:zsolve}
With reference to Notation
\ref{not:a}, for mutually distinct $i,j,k,\ell \in \I$
we have
\begin{eqnarray*}
z_k &=& \frac{\beta_\ell-\beta_i}{\beta_k-\beta_\ell} z_i
+
 \frac{\beta_\ell-\beta_j}{\beta_k-\beta_\ell} z_j,
\\
z_\ell &=& \frac{\beta_i-\beta_k}{\beta_k-\beta_\ell} z_i
+
 \frac{\beta_j-\beta_k}{\beta_k-\beta_\ell} z_j.
\end{eqnarray*}
\end{lemma}

\begin{lemma}
\label{lem:manyb}
With reference to Notation
\ref{not:a}, 
for distinct $i,j \in \I$ the elements
\begin{eqnarray*}
z^r_iz^s_j \qquad \qquad 0 \leq r,s< \infty
\end{eqnarray*}
form a basis for $\mathcal P$.
\end{lemma}

\begin{lemma}
\label{lem:6b}
With reference to Notation
\ref{not:a}, 
for an integer $d \geq 0$ and
distinct $i,j \in \I$ 
 the elements
$\lbrace z^{d-n}_i z^n_j\rbrace_{n=0}^d$ form a basis for
${\mathcal P}_d$.
\end{lemma}

\noindent 
Our next goal is to display a
$\boxtimes$-module
structure on
$\mathcal P$.
The following definition will be useful.

\begin{definition}
\label{def:der}
\rm
\cite[p.~4]{humphreys}
With reference to Notation
\ref{not:a}, by a {\it derivation} of $\mathcal P$
we mean an $\K$-linear map $D :{\mathcal P}\to {\mathcal P}$
such that
$D(uv)= D(u)v+uD(v)$ for $u,v\in {\mathcal P}$.
\end{definition}

\noindent We have a comment.

\begin{lemma}
\label{lem:dercom}
With reference to Notation
\ref{not:a} let $D$ denote a derivation of
$\mathcal P$.
Then $D(z^r)=rz^{r-1}D(z)$ 
for $z\in {\mathcal P}$ and $r\geq 1$.
Moreover $D(1)=0$.
\end{lemma}
\noindent {\it Proof:} The first equation
is a routine consequence of Definition \ref{def:der}.
Setting $z=1$ and $r=2$ in this equation we find
$D(1)=2D(1)$ so $D(1)=0$.
\hfill $\Box$ \\

\begin{lemma}
\label{def:derzero}
With reference to Notation
\ref{not:a}, a derivation of
$\mathcal P$ is 0 if and only if it vanishes on 
${\mathcal P}_1$.
\end{lemma}
\noindent {\it Proof:} 
Let $D$ denote the derivation in question.
We assume $D$ vanishes on ${\mathcal P}_1$ and
show $D=0$.
For distinct $i,j\in \I$ and integers 
$r,s\geq 0$ we have
$D(z^r_i z^s_j)=0$ by
Definition
\ref{def:der} and
Lemma
\ref{lem:dercom}. Now $D=0$ in view of
Lemma
\ref{lem:manyb}.
\hfill $\Box$ \\ 

\begin{proposition}
\label{lem:tetder}
With reference to Notation
\ref{not:a}, there exists a unique $\boxtimes$-module structure
on $\mathcal P$ that satisfies (i), (ii) below:
\begin{enumerate}
\item
each element of $\boxtimes$ acts as a derivation on $\mathcal P$;
\item
$x_{ij}.z_i=-z_i$
and $x_{ij}.z_j=z_j$ for distinct $i,j\in \I$.
\end{enumerate}
\end{proposition}
\noindent {\it Proof:} 
Let $\mbox{Der}\,
{\mathcal P}$
denote the set of derivations for 
${\mathcal P}$. Recall that
 $\mbox{Der}\,
{\mathcal P}$ is a Lie algebra over 
$\K$ with Lie bracket
$\lbrack D,D'\rbrack=DD'-D'D$.
For distinct $i,j \in \I$ define an $\K$-linear
map $D_{ij}:
{\mathcal P}\to 
{\mathcal P}$
such that
\begin{eqnarray*}
D_{ij}(z^r_iz^s_j) =
(s-r)z^r_iz^s_j  
\qquad \qquad r,s=0,1,2,\ldots
\end{eqnarray*}
One checks that $D_{ij}$ is the unique
element of 
 $\mbox{Der}\,
{\mathcal P}$ that sends $z_i\mapsto -z_i$
and $z_j\mapsto z_j$. 
Next one checks that the maps
$\lbrace D_{ij} |i,j\in \I,i\not=j\rbrace $ satisfy the
defining relations for $\boxtimes$ given in
Definition
\ref{def:tet}. To do this, it suffices to
verify that these relations hold on ${\mathcal P}_1$, in view
of Lemma
\ref{def:derzero}. From our comments so far,
there exists a Lie algebra homomorphism
from $\boxtimes$ to 
 $\mbox{Der}\,
{\mathcal P}$ that sends $x_{ij} \mapsto D_{ij}$
for all distinct $i,j\in \I$.
This shows that the required $\boxtimes$-module structure exists.
This $\boxtimes$-module structure is unique in view of
Lemma
\ref{def:derzero}.
\hfill $\Box$ \\

\noindent We emphasize the following.

\begin{lemma}
\label{lem:xijact}
With reference to Notation
\ref{not:a} and
Proposition
\ref{lem:tetder}, 
for distinct $i,j\in \I$ 
and integers $r,s\geq 0$
the element $z^r_i z^s_j$
is an eigenvector for $x_{ij}$ with eigenvalue
$s-r$.
\end{lemma}

\begin{lemma}
\label{lem:adirred}
With reference to Notation
\ref{not:a} the following (i)--(iii) hold.
\begin{enumerate}
\item For $d\geq 0$ the subspace ${\mathcal P}_d$
is an irreducible $\boxtimes$-submodule of $\mathcal P$.
\item 
The $\boxtimes$-module ${\mathcal P}_0$ is trivial.
\item For $d\geq 1$ the $\boxtimes$-module
${\mathcal P}_d$ is isomorphic to an
evaluation module.
\end{enumerate}
\end{lemma}
\noindent {\it Proof:} 
Pick distinct $i,j\in \I$ and consider 
the basis
$\lbrace z^{d-n}_iz^n_j \rbrace_{n=0}^d$ for
${\mathcal P}_d$ from Lemma
\ref{lem:6b}. By 
Lemma
\ref{lem:xijact},
for $0 \leq n \leq d$ the vector
$z^{d-n}_iz^n_j$ is an eigenvector for
$x_{ij}$ with eigenvalue $2n-d$.
Therefore ${\mathcal P}_d$ is
invariant under $x_{ij}$.
We have now shown that
${\mathcal P}_d$ is a $\boxtimes$-submodule.
This module is irreducible
since the eigenvalues of
$x_{ij}$ on ${\mathcal P}_d$ are
$2n-d$ $(0 \leq n \leq d)$
and the corresponding eigenspaces have dimension 1.
This proves (i). To get
(ii) note that ${\mathcal P}_0$ has dimension 1.
To get (iii) apply
Proposition
\ref{prop:evchar3}.
\hfill $\Box$ \\ 

\noindent We now find the evaluation parameter for
the $\boxtimes$-modules given in 
Lemma \ref{lem:adirred}(iii).
The following lemma will
be useful.

\begin{lemma}
\label{lem:derdep}
With reference to Notation
\ref{not:a},
for mutually distinct $i,j,k,\ell \in \I$
the following vanishes on the $\boxtimes$-module $\mathcal P$:
\begin{eqnarray}
\label{eq:betaexp}
(\beta_i-\beta_j)
(\beta_k-\beta_\ell)x_{ij}
+
(\beta_i-\beta_k)
(\beta_\ell-\beta_j)x_{ik}
+
(\beta_i-\beta_\ell)
(\beta_j-\beta_k)x_{i\ell}.
\end{eqnarray}
\end{lemma}
\noindent {\it Proof:} 
Using Lemma
\ref{lem:zsolve} and
Proposition 
\ref{lem:tetder}(ii) we routinely find
that the derivation 
(\ref{eq:betaexp}) vanishes on each of $z_i, z_j$ and
hence on ${\mathcal P}_1$. This derivation is 0
in view of
Lemma
\ref{def:derzero}.
\hfill $\Box$ \\

\begin{theorem}
\label{lem:irred}
With reference to Notation
\ref{not:a}, for 
an integer $d\geq 1$ the 
 $\boxtimes$-module 
${\mathcal P}_d$ is isomorphic to $V_d(a)$ where
\begin{eqnarray*}
a = \frac{\beta_0-\beta_1}{\beta_0-\beta_3}
\, \frac{\beta_2-\beta_3}{\beta_2-\beta_1}.
\end{eqnarray*}
\end{theorem}
\noindent {\it Proof:} 
We invoke
Proposition \ref{prop:evchar} with $V={\mathcal P}_d$.
Note that $a\not\in \lbrace 0,1\rbrace$ since
the $\beta_i$ $(i\in \I)$ are mutually distinct.
By Lemma
\ref{lem:derdep},
the four expressions in
(\ref{eq:ad1}), 
(\ref{eq:ad2}) 
vanish on
$\mathcal P$ and in particular on
${\mathcal P}_d$. By  this and 
Proposition \ref{prop:evchar} the $\boxtimes$-module
${\mathcal P}_d$ 
has evaluation
parameter $a$.  
Now the $\boxtimes$-module
${\mathcal P}_d$ 
is isomorphic to $V_d(a)$
since 
${\mathcal P}_d$ has dimension $d+1$.
\hfill $\Box$ \\

\noindent Earlier in the paper we described the
$\boxtimes$-module $V_d(a)$. We now 
consider how things look from the point of
view of ${\mathcal P}_d$.

\begin{proposition}
\label{prop:adec}
With reference to Notation
\ref{not:a},
for an integer $d\geq 0$ and distinct $i,j\in \I$
the decomposition $\lbrack i,j \rbrack $ on ${\mathcal P}_d$
is described as follows: for $0\leq n\leq d$
the $n$th component is spanned by
$z^{d-n}_iz^n_j$.
\end{proposition}
\noindent {\it Proof:} 
By 
Lemma \ref{lem:6b}
the vectors 
$\lbrace z^{d-n}_iz^n_j\rbrace_{n=0}^d$
form a basis for ${\mathcal P}_d$.
By Lemma 
\ref{lem:xijact}, for $0 \leq n \leq d$ the vector
$z^{d-n}_iz^n_j$ is an eigenvector for
$x_{ij}$ with eigenvalue $2n-d$.
The result follows.
\hfill $\Box$ \\

\begin{proposition}
\label{prop:aflag}
With reference to Notation
\ref{not:a},
for an integer $d\geq 0$ and  $i\in \I$
the flag $\lbrack i \rbrack $ on ${\mathcal P}_d$
is described as follows: for $0\leq n\leq d$
the $n$th component is
$z^{d-n}_i{\mathcal P}_n$.
\end{proposition}
\noindent {\it Proof:} 
Let $U_n$ denote the component in question and
pick $j \in \I, j\not=i$.
By  construction $U_n$ is the sum
of components $0,1,\ldots, n$
for the decomposition
$\lbrack i,j\rbrack$ of ${\mathcal P}_d$.
By this and
Proposition
\ref{prop:adec},
$U_n$ has a basis 
$\lbrace z_i^{d-r}z_j^r\rbrace_{r=0}^n$.
Note that $\lbrace z_i^{n-r}z_j^r \rbrace_{r=0}^n$
is a basis for ${\mathcal P}_n$ so
$\lbrace z_i^{d-r}z_j^r\rbrace_{r=0}^n$ is a basis for
$z_i^{d-n}{\mathcal P}_n$.
Therefore
$U_n=
z_i^{d-n}{\mathcal P}_n$.
\hfill $\Box$ \\

\begin{definition}
\label{def:aeta}
With reference to Notation
\ref{not:a}, for an integer $d\geq 1$
and $i \in \I$ 
the element $z_i^d$ is in component 0 of the
flag $\lbrack i \rbrack$ on ${\mathcal P}_d$.
In 
Notation \ref{not}
we let $\eta_i$ denote any 
nonzero element in this component.
For the rest of this section we choose 
 $\eta_i=z^d_i$.
\end{definition}

\begin{proposition}
\label{prop:brackbasis}
With reference to Notation
\ref{not:a},
for an integer $d\geq 1$ and for mutually distinct
$i,j,k,\ell\in \I$ the basis
$\lbrack i,j,k,\ell\rbrack$ of $\mathcal{P}_d$ is
described  as follows.
For $0 \leq n \leq d$ the $n$th component is
\begin{eqnarray}
\label{eq:un}
z_k^{d-n}z_{\ell}^n {d \choose n} \frac{(\beta_j-\beta_k)^{d-n}
(\beta_j-\beta_{\ell})^n}
{
(\beta_i-\beta_j)^d}.
\end{eqnarray}
\end{proposition}
\noindent {\it Proof:} 
Abbreviate $u_n$ 
for the polynomial (\ref{eq:un}) and note
that
$\lbrace u_n\rbrace_{n=0}^d$ 
is a basis for ${\mathcal P}_d$ by
Lemma
\ref{lem:6b}.
We show that this basis satisfies
 Lemma
\ref{def:eta}(i),(ii).
By Proposition
\ref{prop:adec},
for $0 \leq n \leq d$ the vector
$u_n$ is in
component $n$ of the decomposition
$\lbrack k,\ell\rbrack$ for ${\mathcal P}_d$.
Also 
$\eta_i = z_i^d$ by 
Definition \ref{def:aeta}, 
and
\begin{eqnarray*}
z_i^d &=&
\biggl(
 \frac{\beta_j-\beta_k}{\beta_i-\beta_j} z_k
+
 \frac{\beta_j-\beta_\ell}{\beta_i-\beta_j} z_\ell \biggr)^d
 \\
&=& \sum_{n=0}^d u_n
\end{eqnarray*}
by Lemma
\ref{lem:zsolve} and 
 the binomial theorem.
Therefore the basis $\lbrace u_n\rbrace_{n=0}^d$ satisfies
 Lemma
\ref{def:eta}(i),(ii) and the result follows.
\hfill $\Box$ \\

\noindent We have a comment on the $S_4$ action.

\begin{lemma}
\label{lem:s4rev}
With reference to Notation
\ref{not:a} and
Theorem
\ref{lem:irred},
for mutually distinct
$i,j,k,\ell\in \I$ the 
$(i,j,k,\ell)$-relative of $a$ is
\begin{eqnarray*}
 \frac{\beta_i-\beta_{\ell}}{\beta_i-\beta_k}
\, \frac{\beta_j-\beta_k}{\beta_j-\beta_{\ell}}.
\end{eqnarray*}
\end{lemma}
\noindent {\it Proof:} 
Very similar to the proof of
Lemma
\ref{lem:bar}.
\hfill $\Box$ \\ 

\noindent We now consider the standard bilinear form.

\begin{lemma}
\label{lemma:bilaz}
With reference to Notation
\ref{not:a}, for an integer $d\geq 1$
there exists a standard bilinear
form 
$\langle \,,\,\rangle$ 
on the $\boxtimes$-module
${\mathcal P}_d$ that satisfies the following.
For distinct
 $i,j \in \I$,
\begin{eqnarray*}
\langle z_i^d,z_j^d\rangle = (\beta_k-\beta_{\ell})^d
\end{eqnarray*}
where the set $\lbrace k,\ell \rbrace$ is the complement of
$\lbrace i,j\rbrace$ in  $\I$, and the pair $k,\ell$ 
is ordered so that
the sequence $(i,j,k,\ell)$ is sent to 
$(2,0,1,3)$
by an even permutation in $S_4$.
\end{lemma}
\noindent {\it Proof:} 
For the time being  let $\langle \,,\,\rangle$
denote any standard bilinear form on
${\mathcal P}_d$.
For mutually distinct $i,j,k,\ell \in \I$
define
\begin{eqnarray*}
\tau(i,j,k,\ell)= \frac{\langle z^d_i,z^d_j\rangle}{(\beta_k-\beta_\ell)^d},
\end{eqnarray*}
and note that 
$\tau(i,j,k,\ell)\not=0$ by
Lemma
\ref{lem:fudge}(ii).
We claim
\begin{itemize}
\item
$\tau(j,i,k,\ell) = (-1)^d 
\tau(i,j,k,\ell)$,
\item
$\tau(i,k,j,\ell) = (-1)^d 
\tau(i,j,k,\ell)$,
\item
$\tau(i,j,\ell,k) = (-1)^d 
\tau(i,j,k,\ell)$.
\end{itemize}
The first bullet follows from the symmetry/antisymmetry of
$\langle \,,\,\rangle $ in Lemma
\ref{lem:bil}. The third bullet follows from the construction.
To verify the second bullet we show 
\begin{eqnarray}
\frac{\langle z^d_i,z^d_k\rangle}{(\beta_\ell-\beta_j)^d}
=
\frac{\langle z^d_i,z^d_j\rangle}{(\beta_k-\beta_\ell)^d}.
\label{eq:zz}
\end{eqnarray}
To verify 
(\ref{eq:zz}), in the left-hand side eliminate
$z_k$ using the first equation of
Lemma
\ref{lem:zsolve}, and note that
$\langle z_i^d, z^{d-n}_iz^n_j\rangle = 0$ 
for $0 \leq n \leq d-1$ in view of
Lemma \ref{lem:dual} and
Proposition \ref{prop:adec}.
This proves
(\ref{eq:zz})
and we have now verified the three bullets.
Multiplying $\langle \,,\,\rangle $ by a nonzero 
scalar in $\K$ if necessary, we may assume
$\tau(2,0,1,3)=1$. For $\sigma \in S_4$
define $\mbox{sgn}(\sigma)$ to be $1$ if $\sigma$ is
even and $-1$ if $\sigma$ is odd.
Using the three bullets and
$\tau(2,0,1,3)=1$ we obtain 
\begin{eqnarray*}
\tau(i,j,k,\ell)=\mbox{sgn}(\sigma)^d
\end{eqnarray*}
where $\sigma \in S_4$ sends the sequence
$(i,j,k,\ell)$ to $(2,0,1,3)$.
The result follows.
\hfill $\Box$ \\

\begin{corollary}
Referring to the bilinear form $\langle \,,\,\rangle$
in Lemma  
\ref{lemma:bilaz}, for  $i,j\in \I$
the scalar $\langle z^d_i, z^d_j\rangle $ is displayed in
row $i$, column $j$ of the table below.

\bigskip

\centerline{
\begin{tabular}[t]{c|c c c c}
       & $0$ &$1$ &$2$ &$3$ 
 \\ \hline  
$0$ 
& 
$0$ 
&
$(\beta_2-\beta_3)^d$
&    
$(\beta_3-\beta_1)^d$
&    
 $(\beta_1-\beta_2)^d$     
	\\
$1$
&
$(\beta_3-\beta_2)^d$
&
$0$
&
$(\beta_0-\beta_3)^d$
&    
$(\beta_2-\beta_0)^d$
	\\
$2$
&
$(\beta_1-\beta_3)^d$
&
$(\beta_3-\beta_0)^d$
&
$0$
&
$(\beta_0-\beta_1)^d$
              \\
$3$ &
$(\beta_2-\beta_1)^d$
&
$(\beta_0-\beta_2)^d$
&
$(\beta_1-\beta_0)^d$
&
$0$
	\end{tabular}}

\bigskip
\end{corollary}
\noindent {\it Proof:} 
Routine using 
Lemma 
\ref{lem:fudge}(i) and
Lemma
\ref{lemma:bilaz}.
\hfill $\Box$ \\

\begin{theorem}
\label{prop:abil}
With reference to Notation
\ref{not:a}, for an integer $d\geq 1$
consider the standard bilinear
form 
$\langle \,,\,\rangle$ 
on 
${\mathcal P}_d$ from Lemma
\ref{lemma:bilaz}.
For distinct $i,j \in \I$ and  $0 \leq r,s \leq d$ we have
\begin{eqnarray}
\label{eq:bilinear}
\langle z_i^{d-r}z_j^r,
z_i^{d-s}z_j^s\rangle
= \delta_{r+s,d}(-1)^r {d \choose r}^{-1} (\beta_k-\beta_{\ell})^d
\end{eqnarray}
where the set $\lbrace k,\ell \rbrace$ is the complement of
$\lbrace i,j\rbrace$ in  $\I$, and the pair $k,\ell$ 
is ordered so that the sequence
$(i,j,k,\ell)$ is sent to 
$(2,0,1,3)$ 
by an even permutation in $S_4$.
\end{theorem}
\noindent {\it Proof:} 
Assume $r+s=d$; otherwise
(\ref{eq:bilinear})
holds
by
Lemma \ref{lem:dual} and
Proposition \ref{prop:adec}.
Let $\lbrace u_n\rbrace_{n=0}^d$ denote the basis
$\lbrack k,\ell,i,j\rbrack$ for ${\mathcal P}_d$.
 Using 
 Proposition 
\ref{lem:fins}
 and Proposition
\ref{prop:brackbasis} 
  we find
\begin{eqnarray*}
(-1)^r {d \choose r} &=&
\frac{\langle u_r,u_{d-r}\rangle} 
{\langle u_0,u_d\rangle}
\\
&=& 
\frac{\langle z_i^{d-r}z_j^r,
z_i^{r}z_j^{d-r}\rangle}
{\langle z_i^d,
z_j^d\rangle}
{d \choose r}^2.
\end{eqnarray*}
Line 
(\ref{eq:bilinear})
follows 
in view of Lemma \ref{lemma:bilaz}.
\hfill $\Box$ \\

\noindent Recall the subgroup $G$ of $S_4$
from
Definition \ref{def:G}.
In Corollary 
\ref{thm:n2} we showed that
for an evaluation
module of $\boxtimes$ and for $\sigma \in G$,
twisting the module via $\sigma$
does not change the isomorphism class of the module.
We now interpret this fact using
the $\boxtimes$-module $\mathcal P$.

\medskip
\noindent 
With reference to Notation
\ref{not:a}, recall an {\it automorphism} of ${\mathcal P}$
is an $\K$-linear bijection $\phi: {\mathcal P}\to
{\mathcal P}$ such that
$\phi(uv)=\phi(u)\phi(v)$ for $u,v\in {\mathcal P}$.

\begin{lemma}
\label{lem:compder}
With reference to Notation
\ref{not:a}, for an automorphism $\phi$ 
of $\mathcal P$ and a  derivation $D$ of
$\mathcal P$ 
the composition
$\phi^{-1}D\phi$ is a derivation of 
$\mathcal P$.
\end{lemma}
\noindent {\it Proof:} 
Routine using Definition
\ref{def:der}
and the definition of automorphism.
\hfill $\Box$ \\ 

\begin{lemma}
\label{lem:actsas}
With reference to Notation
\ref{not:a}, 
for mutually distinct $i,j,k,\ell \in \I$ there exists 
a unique automorphism of $\mathcal P$ that sends
\begin{eqnarray}
\label{eq:lineone}
&&z_i \mapsto \frac{\beta_j-\beta_k}{\beta_i-\beta_k}z_j,
\qquad \qquad 
z_j \mapsto \frac{\beta_i-\beta_\ell}{\beta_j-\beta_\ell}z_i,
\\
&&z_k \mapsto \frac{\beta_\ell-\beta_i}{\beta_i-\beta_k}z_\ell,
\qquad \qquad 
z_\ell \mapsto  \frac{\beta_k-\beta_j}{\beta_j-\beta_\ell}z_k.
\label{eq:linetwo}
\end{eqnarray}
\end{lemma}
\noindent {\it Proof:} 
Since $z_i, z_j$ form a basis for ${\mathcal P}_1$
there exists 
an  $\F$-linear transformation
$
\phi :
{\mathcal P}_1 \to
{\mathcal P}_1$
that satisfies
(\ref{eq:lineone}).
Observe that $\phi^{-1}$ exists, and that
$\phi$ extends uniquely to 
an automorphism of $\mathcal P$.
To get 
(\ref{eq:linetwo})
combine 
Lemma \ref{lem:zsolve}
and (\ref{eq:lineone}).
\hfill $\Box$ \\ 

\begin{proposition}
\label{prop:gsig}
With reference to Notation
\ref{not:a}, 
for $\sigma \in G$ there exists an automorphism
$\phi_\sigma$ of $\mathcal P$ 
that sends 
$z_r$ into $\F z _{\sigma(r)}$ for all $r \in \I$.
\end{proposition}
\noindent {\it Proof:} 
Assume $\sigma$ is not the identity; otherwise
the result is clear.
By Definition
\ref{def:G} there exists mutually distinct
$i,j,k,\ell \in \I$ such that
$\sigma=(i,j)(k,\ell)$; let $\phi_\sigma$ denote the
corresponding automorphism of $\mathcal P$ from
Lemma
\ref{lem:actsas}. 
By 
(\ref{eq:lineone}), 
(\ref{eq:linetwo}) the automorphism
$\phi_\sigma$ sends
$z_r$ into $\F z _{\sigma(r)}$ for all $r \in \I$.
\hfill $\Box$ \\ 

\begin{proposition}
\label{prop:similar}
With reference to Notation
\ref{not:a}, 
for $\sigma \in G$ and $\xi\in \boxtimes$
the equation 
\begin{eqnarray}
\label{eq:sim}
\sigma(\xi)=\phi_\sigma \xi \phi^{-1}_\sigma
\end{eqnarray}
holds on $\mathcal P$.
\end{proposition}
\noindent {\it Proof:} 
Without loss we may assume that $\xi$ is a generator
$x_{rs}$ of $\boxtimes$.
Each side of (\ref{eq:sim})
acts on $\mathcal P$ as a derivation,
by
Proposition
\ref{lem:tetder} and
Lemma
\ref{lem:compder}. 
Now by Lemma
\ref{def:derzero},
it suffices
to show that
 (\ref{eq:sim}) holds on ${\mathcal P}_1$.
The elements $z_{\sigma(r)}$ and
$z_{\sigma(s)}$ form a basis for ${\mathcal P}_1$.
We now apply each side of
 (\ref{eq:sim}) to
$z_{\sigma(r)}$.
Concerning the left-hand side,
\begin{eqnarray*}
\sigma(x_{rs}).z_{\sigma(r)} &=& 
x_{\sigma(r),\sigma(s)}.z_{\sigma(r)}
\\
&=& -z_{\sigma(r)}.
\end{eqnarray*}
Concerning the right-hand side,
first note by Proposition
\ref{prop:gsig} that there exists 
$\gamma \in \K$ such that $\phi_\sigma(z_r) = \gamma z_{\sigma(r)}$.
Observe $\gamma\not=0$ since $\phi_\sigma$ is a bijection.
We have
\begin{eqnarray*}
\phi_\sigma x_{rs} \phi^{-1}_\sigma.z_{\sigma(r)}
&=&
\gamma^{-1} \phi_\sigma x_{rs}. z_r
\\
&=&
-\gamma^{-1} \phi_\sigma. z_r
\\
&=&
- z_{\sigma(r)}.
\end{eqnarray*}
Thus the two sides
 of (\ref{eq:sim}) agree at
 $z_{\sigma(r)}$.
A similar argument shows that
the two sides of
 (\ref{eq:sim}) agree at 
 $z_{\sigma(s)}$. Now 
 (\ref{eq:sim}) holds on
${\mathcal P}_1$ and the result follows. 
\hfill $\Box$ \\

\begin{theorem}
\label{thm:mainth}
With reference to
Notation
\ref{not:a}
and Proposition
\ref{prop:gsig},
 for $\sigma \in G$ the
map $\phi_\sigma$
is an isomorphism of $\boxtimes$-modules
from $\mathcal P$ to $\mathcal P$ twisted via $\sigma$.
\end{theorem}
\noindent {\it Proof:} 
This is a reformulation of
Proposition \ref{prop:similar}.
\hfill $\Box$ \\

\section{Finite-dimensional irreducible $\boxtimes$-modules}

\noindent We have now completed our description of the
evaluation modules for $\boxtimes$.
For the rest of this paper we consider
general finite-dimensional irreducible
$\boxtimes$-modules. For these we will not 
go into as much detail
as we did for the evaluation modules,
but there are a few points we would like to make.

\medskip
\noindent 
We start with some comments on tensor products.

\medskip
\noindent Let $U$ and $V$ denote $\boxtimes$-modules.
By \cite[p.~26]{humphreys}
the vector space $U\otimes V$ has a $\boxtimes$-module
structure given by
\begin{eqnarray}
\label{eq:tens}
\xi.(u\otimes v) = (\xi.u)\otimes v + u\otimes (\xi.v)
\qquad \qquad u \in U,\quad v\in V, \quad \xi \in \boxtimes.
\end{eqnarray}

\begin{lemma}
\label{lem:tcom}
The following (i), (ii) hold for $\boxtimes$-modules
 $U$ and $V$:
\begin{enumerate}
\item There exists an isomorphism of $\boxtimes$-modules
from $U\otimes V$ to $V\otimes U$ that sends
$u \otimes v \mapsto v\otimes u$ for all $u \in U$
and $v \in V$.
\item Assume the 
$\boxtimes$-module $U\otimes V$ is irreducible. Then
$U$ and  $V$ are irreducible.
\end{enumerate}
\end{lemma}
\noindent {\it Proof:} 
(i) Routine.
\\
\noindent (ii) The $\boxtimes$-module $U$ is irreducible
since if $U'$ is a nonzero proper submodule of $U$
then $U'\otimes V$ is a nonzero proper submodule of $U\otimes V$.
The proof for $V$ is similar.
\hfill $\Box $ \\

\begin{theorem}
\label{thm:modclass}
{\rm  
\cite[Section 1]{Ha}}
Every nontrivial finite-dimensional 
irreducible $\boxtimes$-module is isomorphic to a tensor product of
evaluation modules. Two such tensor products are isomorphic
if and only if one can be obtained from the other by
permuting the factors in the tensor product. A tensor product
of evaluation modules
\begin{eqnarray*}
V_{d_1}(a_1) \otimes
V_{d_2}(a_2) \otimes \cdots \otimes
V_{d_N}(a_N)
\end{eqnarray*}
is irreducible if and only if $a_1, a_2, \ldots, a_N$ are mutually
distinct.
\end{theorem}
\noindent {\it Proof:} 
In \cite[Section~1]{Ha} Hartwig
considers a certain Lie algebra $\mathcal O$ 
called the {\it Onsager algebra}
\cite{Onsager}, \cite{perk},
\cite{uglov1}. 
 In
\cite[Theorems~1.7,~1.8]{Ha} he gives an explicit
bijection between
(i) the set of isomorphism classes
of finite-dimensional irreducible
$\boxtimes$-modules;
(ii) the set of
isomorphism classes of finite-dimensional irreducible
$\mathcal O$-modules that have type $(0,0)$.
In \cite[Theorems~1.3,~1.4,~1.6]{Ha} Hartwig
summarizes the classification of 
 finite-dimensional irreducible
$\mathcal O$-modules that have type $(0,0)$. This classification
is due to Davies
\cite{Davfirst},
\cite{Da}; 
see also Date and Roan
\cite{DateRoan2}.
If we interpret this classification as a classification of
the finite-dimensional irreducible $\boxtimes$-modules
via the above bijection, we get the present theorem.
To aid in this interpretation we note that our
definition of the evaluation parameter does not
match Hartwig's definition of the evaluation
parameter. Denoting our evaluation parameter by $a$ and Hartwig's
evaluation parameter by $b$ 
we have $a=4b(b+1)^{-2}$.
\hfill $\Box $ \\

\noindent We have some comments.

\begin{lemma}
\label{lem:tenstwist}
Let $U$ and $V$ denote $\boxtimes$-modules. Then 
for $\sigma \in S_4$ the following coincide:
\begin{enumerate}
\item the $\boxtimes$-module ${}^\sigma (U\otimes V)$;
\item the $\boxtimes$-module ${}^\sigma U\otimes {}^\sigma V$.
\end{enumerate}
\end{lemma}
\noindent {\it Proof:} 
Routine using  Definition
\ref{def:twist}
and 
(\ref{eq:tens}).
\hfill $\Box $ \\

\begin{theorem}
\label{lem:twistdec}
Let $V$ denote a nontrivial finite-dimensional irreducible
$\boxtimes$-module, and write $V$ as a tensor product
of evaluation modules:
\begin{eqnarray*}
V=V_{d_1}(a_1) \otimes
V_{d_2}(a_2) \otimes \cdots \otimes
V_{d_N}(a_N).
\end{eqnarray*}
Then for $\sigma \in S_4$ the $\boxtimes$-module 
${}^\sigma V$ is isomorphic to
\begin{eqnarray*}
V_{d_1}(\sigma(a_1)) \otimes
V_{d_2}(\sigma(a_2)) \otimes \cdots \otimes
V_{d_N}(\sigma(a_N)).
\end{eqnarray*}
\end{theorem}
\noindent {\it Proof:} 
Combine Theorem
\ref{thm:evtwist}
and
Lemma \ref{lem:tenstwist}.
\hfill $\Box $ \\

\begin{corollary}
Let $V$ denote a finite-dimensional irreducible
$\boxtimes$-module and let the group $G$ be as in
Definition 
\ref{def:G}.
 Then for $\sigma \in G$ the following
are isomorphic:
\begin{enumerate}
\item the $\boxtimes$-module $V$ twisted via $\sigma$;
\item the $\boxtimes$-module $V$.
\end{enumerate}
\end{corollary}
\noindent {\it Proof:} 
Combine 
Lemma \ref{thm:n1}
and Theorem 
\ref{lem:twistdec}.
\hfill $\Box $ \\

\section{Finite-dimensional irreducible $\boxtimes$-modules; the
shape polynomial}

\noindent In this section we obtain an explicit formula
for the shape of a 
finite-dimensional irreducible $\boxtimes$-module.
We use the following notation.
For an indeterminate $\lambda$ 
let $\K\lbrack \lambda \rbrack$
denote the $\K$-algebra consisting of all polynomials 
in $\lambda $ that have coefficients in $\K$.

\begin{definition} 
\label{def:sgen}
\rm
Let $V$ denote a finite-dimensional irreducible $\boxtimes$-module.
We define a polynomial $S_V \in \F\lbrack \lambda \rbrack$ by 
\begin{eqnarray*}
S_V = \sum_{n=0}^d \rho_n \lambda^n,
\end{eqnarray*}
where $\lbrace \rho_n\rbrace_{n=0}^d$ is the shape of
$V$. We call $S_V$ the {\it shape polynomial} for $V$.
\end{definition}

\begin{example}
\label{ex:shape}
Let $V$ denote an evaluation module for $\boxtimes$. Then
\begin{eqnarray*}
S_V = 1 + \lambda + \lambda^2 + \cdots+ \lambda^d
\end{eqnarray*}
where $d$ is the diameter of $V$.
\end{example}
\noindent {\it Proof:} 
Combine
Proposition \ref{prop:evchar3}
and Definition
\ref{def:sgen}.
\hfill $\Box $ \\

\noindent
Let $U, V$ denote finite-dimensional irreducible
$\boxtimes$-modules such that the $\boxtimes$-module
$U\otimes V$ is irreducible.
We are going to show that $S_{U\otimes V}=S_US_V$. 
We will use the following lemma.

\begin{lemma}
\label{lem:shapecomb}
Let $U, V$ denote finite-dimensional irreducible
$\boxtimes$-modules such that the $\boxtimes$-module
$U\otimes V$ is irreducible. Let $d$ (resp. $\delta$)
denote the diameter of $U$ (resp. $V$).
\begin{enumerate}
\item The diameter of $U\otimes V$ is $d+\delta$.
\item For distinct $i,j\in \I$ the 
decomposition $\lbrack i,j\rbrack $
of $U\otimes V$ is described as follows. For
 $0 \leq n \leq d+\delta$ 
the $n$th component is
\begin{eqnarray}
\label{eq:uv}
\sum_{r,s} U_r \otimes V_s,
\end{eqnarray}
where 
$\lbrace U_r\rbrace_{r=0}^d$ 
(resp. 
$\lbrace V_s\rbrace_{s=0}^{\delta}$)
denotes the decomposition
$\lbrack i,j\rbrack$ of $U$ (resp. $V$), and
where the sum is over all ordered pairs $r,s$ such
that $0\leq r\leq d$, $0 \leq s \leq \delta$, $r+s=n$.
\end{enumerate}
\end{lemma}
\noindent {\it Proof:} 
For $0 \leq n \leq d+\delta$ let
$(U\otimes V)_n$ denote the sum 
(\ref{eq:uv}). By  
(\ref{eq:tens}), each element of 
$(U\otimes V)_n$ is an eigenvector for 
$x_{ij}$ with eigenvalue $2n-d-\delta$.
By construction the sequence
$\lbrace (U\otimes V)_n\rbrace_{n=0}^{d+\delta}$
is a decomposition of $U\otimes V$.
The results follow.
\hfill $\Box $ \\

\begin{corollary}
\label{cor:shape}
Let $U, V$ denote finite-dimensional irreducible
$\boxtimes$-modules such that the $\boxtimes$-module
$U\otimes V$ is irreducible.
Then with reference to
Definition 
\ref{def:sgen}, 
\begin{eqnarray*}
S_{U\otimes V} = S_U S_V.
\end{eqnarray*}
\end{corollary}
\noindent {\it Proof:} 
Adopt the notation of
Lemma
\ref{lem:shapecomb}.
For $0 \leq n \leq d+\delta$
the sum
(\ref{eq:uv}) is direct so
\begin{eqnarray}
\label{eq:32dir}
\rho_n(U\otimes V) = 
\sum_{r,s} \rho_r(U) \rho_s(V),
\end{eqnarray}
where the sum is over all ordered pairs $r,s$ such
that $0\leq r\leq d$, $0 \leq s \leq \delta$, $r+s=n$.
By Definition
\ref{def:sgen},
\begin{eqnarray*}
S_{U\otimes V} = 
\sum_{n=0}^{d+\delta}
\rho_n(U\otimes V) \lambda^n.
\end{eqnarray*}
Evaluating this using
(\ref{eq:32dir}) we routinely find
$S_{U\otimes V} =S_US_V$. 
\hfill $\Box $ \\

\begin{theorem}
\label{lem:shape}
Let $V$ denote a nontrivial finite-dimensional irreducible
$\boxtimes$-module, and write $V$ as a tensor product
of evaluation modules:
\begin{eqnarray*}
V=V_{d_1}(a_1) \otimes
V_{d_2}(a_2) \otimes \cdots \otimes
V_{d_N}(a_N).
\end{eqnarray*}
Then the shape $\lbrace \rho_n\rbrace_{n=0}^d$ of $V$ satisfies
\begin{eqnarray*}
\sum_{n=0}^d \rho_n \lambda^n = \prod_{j=1}^N 
(1+\lambda+\lambda^2+\cdots + \lambda^{d_j}).
\end{eqnarray*}
In particular
$\rho_0=1$.
\end{theorem}
\noindent {\it Proof:} 
Combine Example
\ref{ex:shape}
and Corollary \ref{cor:shape}.
\hfill $\Box $ \\

\section{Finite-dimensional irreducible $\boxtimes$-modules; the
Drinfel'd polynomial}

\noindent Let $V$ denote a finite-dimensional irreducible
$\boxtimes$-module.  
In this section we define
a certain polynomial $P_V$ called the Drinfel'd polynomial.
We show that the map $V\mapsto P_V$ induces a bijection between
the following two sets: (i)
the isomorphism classes of finite-dimensional irreducible
$\boxtimes$-modules;
(ii) the polynomials in $\K\lbrack \lambda \rbrack$ that
have constant coefficient 1 and are nonzero at $\lambda=1$.
We begin with a comment.

\begin{lemma}
\label{lem:raise}
Let $V$ denote a finite-dimensional irreducible
$\boxtimes$-module. Pick mutually distinct $i,j,k\in\I$ 
and consider the action of $x_{ij}+x_{jk}$ on the 
decomposition $\lbrack i,j\rbrack$ of $V$. Denoting this
decomposition by $\lbrace V_n\rbrace_{n=0}^d$ 
we have
$(x_{ij}+x_{jk})V_n\subseteq V_{n+1}$ for $0 \leq n \leq d$.
\end{lemma}
\noindent {\it Proof:} 
By the definition of 
$\lbrack i,j\rbrack$,
\begin{eqnarray}
\label{eq:ij1}
(x_{ij}-(2n-d)I)V_n=0.
\end{eqnarray}
By the table in Section 3,
\begin{eqnarray}
\label{eq:ij2}
(x_{jk}-(d-2n)I)V_n \subseteq V_{n+1}.
\end{eqnarray}
Adding 
(\ref{eq:ij1}), 
(\ref{eq:ij2})  we get the result.
\hfill $\Box $ \\

\begin{definition}
\label{def:si}
\rm
Let $V$ denote a finite-dimensional irreducible
$\boxtimes$-module and let $\lbrace V_n\rbrace_{n=0}^d$
denote the decompostion $\lbrack 1,3\rbrack$ of $V$.
Note that 
$\mbox{dim}(V_0)=1$ by
Theorem
\ref{lem:shape}.
Abbreviate 
\begin{eqnarray}
\label{eq:epm}
e^+:=\frac{x_{13}+x_{30}}{2},
\qquad \qquad
e^-:=\frac{x_{31}+x_{12}}{2}.
\end{eqnarray}
By Lemma
\ref{lem:raise} and since the decomposition
$\lbrack 3,1\rbrack$ is the
inversion of $\lbrack 1,3\rbrack$,
\begin{eqnarray}
\label{eq:updown}
e^+V_n \subseteq V_{n+1},
\qquad \qquad
e^-V_n \subseteq V_{n-1}
\qquad \qquad (0 \leq n \leq d).
\end{eqnarray}
For an integer $i\geq 0$ the space
$V_0$ is invariant under $(e^-)^i(e^+)^i$; let
$\vartheta_i=\vartheta_i(V)$ denote the corresponding eigenvalue.
Note that 
$\vartheta_i=0$ for $i>d$.
\end{definition}

\begin{definition}
\label{def:dri}
\rm
Let $V$ denote a finite-dimensional irreducible
$\boxtimes$-module. We define a polynomial
$P_V \in \F\lbrack \lambda \rbrack$ by
\begin{eqnarray}
\label{eq:drin}
P_V = \sum_{i=0}^{\infty}
\frac{(-1)^i \vartheta_i \lambda^i}{(i!)^2},
\end{eqnarray}
where the scalars $\vartheta_i$ are from Definition
\ref{def:si}. Observe that $P_V$ has degree at
 most the diameter of $V$. Moreover $P_V$ has
constant coefficient $\vartheta_0=1$.
Following
\cite[Section 3.4]{charp},
\cite[Definition 4.2]{NN},
\cite[Section 4]{degxz}
we call $P_V$ the
{\it Drinfel'd polynomial} of $V$.
\end{definition}

\begin{example}
\label{ex:drin}
Let $V=V_d(a)$ denote an evaluation module
for $\boxtimes$. Then 
$P_V = (1-a\lambda)^d$.
\end{example}
\noindent {\it Proof:} 
Let $\lbrace V_n\rbrace_{n=0}^d$ denote the decomposition
$\lbrack 1,3\rbrack$ of $V$. Let $\lbrace v_n\rbrace_{n=0}^d$
denote the basis for the $\mathfrak{sl}_2$-module
$V_d$ from Lemma
\ref{lem:sl2mod}. Note that 
 $v_n$ is a basis of $V_n$ for $0 \leq n \leq d$.
Using 
(\ref{eq:xyz}),
(\ref{eq:epm})
and Lemma
\ref{lem:evchar} we find
$EV_a(e^{+})=af$ and 
$EV_a(e^{-})=e$.
By this and 
Lemma
\ref{lem:sl2mod} we find
$e^+.v_n=a(n+1)v_{n+1}$ for $0 \leq n \leq d-1$, $e^+.v_d=0$,
$e^-.v_n=(d-n+1)v_{n-1}$ for $1 \leq n \leq d$, $e^-.v_0=0$.
Using this data we find
$\vartheta_i=a^i(i!)^2 {d \choose i}$ for $0 \leq i \leq d$.
Now (\ref{eq:drin}) becomes
\begin{eqnarray*}
P_V &=& 
\sum_{i=0}^d
{d \choose i}(-1)^i a^i \lambda^i
\\
&=& (1-a\lambda)^d
\end{eqnarray*}
by the binomial theorem.
\hfill $\Box$ \\

\noindent The Drinfel'd polynomial has the following
property.

\begin{proposition}
\label{prop:drin}
Let $U,V$ denote finite-dimensional irreducible
$\boxtimes$-modules such that the $\boxtimes$-module
$U\otimes V$ is irreducible.
Then 
\begin{eqnarray*}
P_{U\otimes V} = P_U P_V.
\end{eqnarray*}
\end{proposition}
\noindent {\it Proof:} 
We claim that for an integer $i\geq 0$,
\begin{eqnarray}
\label{eq:sigd}
\vartheta_i(U\otimes V) = \sum_{n=0}^i 
{i \choose n}^2
\vartheta_{i-n}(U)\vartheta_n(V).
\end{eqnarray}
To prove the claim,
let 
$U_0$ (resp. $V_0$) denote the
$0$th component for the decomposition
 $\lbrack 1,3\rbrack$ of $U$ (resp. $V$).
By Theorem \ref{lem:shape},
 each of $U_0$, $V_0$ has dimension 1.
By 
Lemma \ref{lem:shapecomb}(ii) the space
$U_0\otimes V_0$ is
the $0$th component of the decomposition
 $\lbrack 1,3\rbrack$ of $U\otimes V$.
Therefore by Definition
\ref{def:si} the scalar $\vartheta_i(U\otimes V)$ is
the eigenvalue of 
 $(e^-)^i(e^+)^i$ associated with  
$U_0\otimes V_0$.
Pick $ 0 \not=u \in U_0$ and $0\not=v \in V_0$.
Using (\ref{eq:tens}) we obtain
\begin{eqnarray}
(e^-)^i(e^+)^i.(u\otimes v) &=&
(e^-)^i \sum_{n=0}^i {i \choose n} ((e^+)^{i-n}.u)\otimes ((e^+)^n.v)
\nonumber
\\
&=& \sum_{m=0}^i \sum_{n=0}^i
{i \choose m}{i \choose n} ((e^-)^{i-m}(e^+)^{i-n}.u)\otimes
((e^-)^m(e^+)^n.v).
\label{eq:doublesum}
\end{eqnarray}
We examine the terms in 
(\ref{eq:doublesum}).
By (\ref{eq:updown}) and the line below it,
for $0 \leq m,n\leq i$ the vector
$(e^-)^m(e^+)^n.v$ is equal to 0 if 
$m>n$ and $\vartheta_n(V)v$ if $m=n$.
Similarly 
$(e^-)^{i-m}(e^+)^{i-n}.u$ is equal to 0 if
$m<n$ and 
$\vartheta_{i-n}(U)u$ if $m=n$.
By these comments the double sum  in
(\ref{eq:doublesum}) 
is equal to $u\otimes v$ times the sum on the right
in 
(\ref{eq:sigd}). We conclude
that
(\ref{eq:sigd}) is valid and the claim is proved.
By Definition
\ref{def:dri},
\begin{eqnarray*}
P_{U\otimes V} = \sum_{i=0}^\infty
\frac{(-1)^i \vartheta_i(U\otimes V) \lambda^i}{(i!)^2}.
\end{eqnarray*}
Evaluating this using
(\ref{eq:sigd}) we get
$P_{U\otimes V} = P_UP_V$ after a brief calculation.
\hfill $\Box$ \\

\begin{theorem} 
\label{thm:dr}
Let $V$ denote 
a nontrivial
finite-dimensional irreducible $\boxtimes$-module, and
write $V$ as a tensor product of evaluation modules:
\begin{eqnarray*}
V=V_{d_1}(a_1) \otimes
V_{d_2}(a_2) \otimes \cdots \otimes
V_{d_N}(a_N).
\end{eqnarray*}
Then the Drinfel'd polynomial $P_V$ is given by 
\begin{eqnarray*}
P_V = 
\prod_{j=1}^N 
(1-a_j\lambda)^{d_j}.
\end{eqnarray*}
\end{theorem}
\noindent {\it Proof:} 
Combine Example
\ref{ex:drin} and 
Proposition
\ref{prop:drin}.
\hfill $\Box $ \\

\begin{corollary}
The map $V \mapsto P_V$ induces a bijection between the following
two sets:
\begin{enumerate}
\item the isomorphism classes of finite-dimensional irreducible
$\boxtimes$-modules;
\item the polynomials in $\K\lbrack \lambda \rbrack$ that
have constant coefficient 1 and are nonzero at $\lambda=1$.
\end{enumerate}
\end{corollary}
\noindent {\it Proof:} 
Combine
Theorem
\ref{thm:modclass}
and
Theorem \ref{thm:dr}.
\hfill $\Box$ \\

\section{Finite-dimensional irreducible $\boxtimes$-modules; the
dual space}

\noindent In this section we show that 
a finite-dimensional irreducible
$\boxtimes$-module is isomorphic to its dual.

\begin{theorem}
\label{thm:form1}
Let $V$ denote a 
finite-dimensional irreducible $\boxtimes$-module.
Then there exists a nonzero
bilinear form $\langle\,,\,\rangle$ on $V$ such that
\begin{eqnarray}
\langle \xi.u,v\rangle = -\langle u,\xi.v\rangle
\qquad \qquad \xi \in \boxtimes, \qquad u,v\in V.
\label{eq:fullbil}
\end{eqnarray}
This form is unique up to multiplication by a nonzero scalar in 
$\F$.
The form is nondegenerate.
The form is symmetric (resp. antisymmetric)
when the diameter is even (resp. odd).
\end{theorem}
\noindent {\it Proof:} 
Assume $V$ is nontrivial; otherwise the proof is routine.
Concerning the existence of $\langle \,,\,\rangle$
 write $V$ as a tensor product
of evaluation modules:
\begin{eqnarray*}
V = 
V_{d_1}(a_1)\otimes
V_{d_2}(a_2)\otimes \cdots \otimes
V_{d_N}(a_N).
\end{eqnarray*}
For $1 \leq i \leq N$ let
$\langle \,,\,\rangle_i$ denote a standard bilinear
form on 
$V_{d_i}(a_i)$.
Define a bilinear
form
$\langle \,,\,\rangle$ on $V$ such that
\begin{eqnarray}
\langle 
\otimes_{i=1}^N u_i,
\otimes_{i=1}^N v_i\rangle
=
\prod_{i=1}^N \langle u_i, v_i \rangle_i
\label{eq:proddef}
\end{eqnarray}
for all $u_i, v_i \in V_{d_i}(a_i)$  $(1 \leq i \leq N)$.
By construction 
$\langle \,,\,\rangle$  is nonzero.
Using 
(\ref{eq:tens}) 
and
(\ref{eq:proddef}) one checks that
$\langle \,,\,\rangle$
satisfies
(\ref{eq:fullbil}).
To show that
$\langle \,,\,\rangle$ is nondegenerate,
note that the subspace 
$\lbrace u \in V|
\langle u,v\rangle=0\; \forall v \in V\rbrace $
is a proper $\boxtimes$-submodule of $V$
and therefore zero by the irreducibility of $V$.
Concerning the uniqueness of
$\langle \,,\,\rangle$,
let $\langle \,,\,\rangle'$
denote any bilinear form
on $V$ that satisfies
(\ref{eq:fullbil}). 
We show that
$\langle \,,\,\rangle'$ is a scalar multiple
of 
$\langle \,,\,\rangle$.
Pick a basis for $V$,
and for $\xi \in \boxtimes$ let $\xi_b$ denote the 
matrix that represents $\xi$ with respect to
this basis.
Let $M$ (resp. $N$) denote the matrix that represents 
$\langle \,,\,\rangle$  
(resp. 
$\langle \,,\,\rangle'$)  
with respect to the basis.
Note that $M$ is invertible since 
$\langle \,,\,\rangle$  is nondegenerate.
By 
(\ref{eq:fullbil}) we have 
$\xi^t_b M = -M\xi_b$
and  $\xi^t_b N = -N\xi_b$ for $\xi \in \boxtimes$.
Combining these equations we find that
$M^{-1}N$ 
commutes with $\xi_b$ for all $\xi \in \boxtimes$.
Now $M^{-1}N$ is a scalar multiple of the identity
by Schur's lemma 
\cite[Lemma 27.3]{CR}
and since the
$\boxtimes$-module $V$ is irreducible.
By these comments $N$ is a scalar multiple of $M$
so
$\langle \,,\,\rangle'$  is a scalar multiple of
$\langle \,,\,\rangle$.
We have now shown that
$\langle \,,\,\rangle$ is unique up to multiplication
by a nonzero scalar in $\K$.
To verify our assertions concerning
symmetry/asymmetry, let $d$ denote the
diameter of $V$ and note that
$d=\sum_{i=1}^N d_i$ by
Theorem \ref{lem:shape}.
Referring to
(\ref{eq:proddef}),
\begin{eqnarray*}
\langle \otimes_{i=1}^N u_i,
 \otimes_{i=1}^N v_i\rangle
&=& 
\prod_{i=1}^N \langle u_i, v_i \rangle_i
\\
&=& 
\prod_{i=1}^N (-1)^{d_i}\langle v_i, u_i \rangle_i 
\qquad \qquad (\mbox{by Lemma  
\ref{lem:bil}})
\\
&=&(-1)^d\prod_{i=1}^N \langle v_i, u_i \rangle_i
\\
&=& (-1)^d\langle \otimes_{i=1}^N v_i,
 \otimes_{i=1}^N u_i\rangle.
\end{eqnarray*}
It follows that
$\langle \,,\,\rangle$ is symmetric (resp. antisymmetric)
when $d$ is even (resp. $d$ is odd).
\hfill $\Box$ \\ 

\noindent 
Before we state our next result we recall 
a concept.
Let $V$ denote a finite-dimensional vector
space over $\F$.
By definition the dual space
$V^*$ is the vector space over $\K$ consisting
of the linear transformations from $V$
to $\K$. The dimensions of $V$ and $V^*$ coincide.
Now assume that $V$ supports a $\boxtimes$-module structure.
Then $V^*$ carries a $\boxtimes$-module structure
such that for $\xi \in \boxtimes$ and $f \in V^*$,
\begin{eqnarray}
\label{eq:dualprop}
(\xi.f)(v) = -f(\xi.v)      \qquad \qquad v \in V.
\end{eqnarray}

\begin{theorem}
Let $V$ denote a finite-dimensional 
irreducible $\boxtimes$-module and 
let $\langle \,,\,\rangle$ denote
a bilinear form on $V$ from 
Theorem
\ref{thm:form1}.
 Then
there exists an isomorphism of 
$\boxtimes$-modules $\varphi:V \to V^*$
such that
\begin{eqnarray}
\label{eq:gamprop}
\varphi(u)(v) = \langle u,v\rangle 
\qquad \qquad u,v\in V.
\end{eqnarray}
\end{theorem}
\noindent {\it Proof:} 
By elementary linear algebra
 there exists a unique
linear transformation $\varphi :V \to V^*$
that satisfies
(\ref{eq:gamprop}). 
The kernel of $\varphi$ is
$\lbrace u \in V|
\langle u,v\rangle=0\; \forall v \in V\rbrace $.
This space is zero since
$\langle \,,\,\rangle$ is nondegenerate, so
$\varphi$ is injective.
Now since $V,V^*$ have the
same dimension, 
the map $\varphi$ is a bijection and
hence an isomorphism of vector
spaces.
Using
(\ref{eq:fullbil}) and 
(\ref{eq:dualprop}) one checks
that $\varphi$ is an isomorphism
of $\boxtimes$-modules.
\hfill $\Box$ \\ 

\section{Finite-dimensional irreducible $\boxtimes$-modules; 
more  bilinear forms}

Let $V$ denote a finite-dimensional irreducible
$\boxtimes$-module. In Theorem
\ref{thm:form1} we displayed a  bilinear
form $\langle \,,\,\rangle$ on $V$. In this section we mention some
related bilinear forms that are of interest.

\begin{theorem}
\label{thm:bilex}
Let $V$ denote a finite-dimensional irreducible
$\boxtimes$-module.
For a nonidentity $\sigma \in G$ there exists
a nonzero bilinear form $\langle \,,\,\rangle_{\sigma}$
on $V$ such that
\begin{eqnarray}
\label{eq:bilsig}
\langle \xi.u,v\rangle_\sigma = 
- 
\langle u,\sigma(\xi).v\rangle_\sigma
\qquad \qquad 
\xi \in \boxtimes, \qquad u,v \in V.
\end{eqnarray}
This form is unique up to
 multiplication
by a nonzero scalar in $\K$.
The form is nondegenerate.
The form is symmetric.
\end{theorem}
\noindent {\it Proof:} 
Concerning existence, let 
$\langle \,,\,\rangle$ denote
a bilinear form on
$V$ from
Theorem
\ref{thm:form1}.
Let $\zeta :V\to V$ denote an
isomorphism of $\boxtimes$-modules from
$V$ to $V$ twisted via $\sigma$.
Define a bilinear form 
$\langle \,,\,\rangle_\sigma$ on $V$ such that 
$\langle u,v\rangle_\sigma = 
\langle u,\zeta(v)\rangle$ for
all $u,v\in V$.
One checks that 
$\langle \,,\,\rangle_\sigma $
is nonzero and satisfies
(\ref{eq:bilsig}).
Concerning uniqueness, let 
$\langle \,,\,\rangle'_\sigma $ denote
a bilinear form on $V$ that
satisfies 
(\ref{eq:bilsig}). We show that 
$\langle \,,\,\rangle'_\sigma $ is a scalar
multiple of 
$\langle \,,\,\rangle_\sigma$.
Define a bilinear form
$\langle \,,\,\rangle'$ on $V$ by
$\langle u,v\rangle' = 
\langle u,\zeta^{-1}(v)\rangle'_\sigma$ for all
$u,v\in V$. Then 
$\langle \,,\,\rangle'$
satisfies
(\ref{eq:fullbil}). Now by the uniqueness of
$\langle \,,\,\rangle$ there exists
$\beta \in \F$ such that
$\langle \,,\,\rangle'
=\beta \langle \,,\,\rangle$.
This implies that 
$\langle \,,\,\rangle'_\sigma=\beta
\langle \,,\,\rangle_\sigma$.
The form  
$\langle \,,\,\rangle_\sigma$ is nondegenerate
since the subspace 
$\lbrace u \in V|
\langle u,v\rangle_\sigma= 0 \;\forall v \in V\rbrace$
is a proper $\boxtimes$-submodule of $V$ and therefore zero
by the irreducibility of $V$.
We now show that
$\langle \,,\,\rangle_\sigma$ is symmetric.
Define a bilinear form
$\langle \,,\,\rangle^\dagger_\sigma$
on $V$ by
$\langle u,v\rangle^\dagger_\sigma = 
\langle v,u\rangle_\sigma$ for all $u,v\in V$.
By construction and since $\sigma^2=1$ we find
$\langle \,,\,\rangle^\dagger_\sigma$ satisfies
(\ref{eq:bilsig}). Now by the uniqueness of
$\langle \,,\,\rangle_\sigma$ 
there exists
$\gamma \in \F$ such that
$\langle \,,\,\rangle^\dagger_
\sigma=
\gamma \langle \,,\,\rangle_\sigma$.
In other words
$\langle v,u\rangle_\sigma= 
\gamma \langle u,v\rangle_\sigma$
for all $u,v\in V$.
We  show $\gamma=1$.
By Definition
\ref{def:G} there exist mutually distinct
$i,j,k,\ell \in \I$ such that
$\sigma=(i,j)(k,\ell)$.
Let $\lbrace V_n\rbrace_{n=0}^d $ denote the
decomposition $\lbrack i,j\rbrack$ of $V$,
and recall that $V_0$ has dimension 1 by
Theorem
\ref{lem:shape}.
Setting $\xi=x_{ij}$ in
(\ref{eq:bilsig}) and using $\sigma(x_{ij})=x_{ji}=-x_{ij}$
we find
$\langle x_{ij}.u,v\rangle_\sigma
=
\langle u,x_{ij}.v\rangle_\sigma $
for all $u,v\in V$.
It follows that
 $\lbrace V_n\rbrace_{n=0}^d $ 
are mutually orthogonal
with respect to
$\langle \,,\,\rangle_\sigma$.
By this and since 
$\langle \,,\,\rangle_\sigma$ is nondegenerate we find
the restriction of
$\langle \,,\,\rangle_\sigma$ to $V_0$ is nonzero.
By our above comments,
for nonzero $u \in V_0$ we have
$\langle u,u\rangle_\sigma \not=0$
 and $\langle u,u\rangle_\sigma 
= \gamma 
\langle u,u\rangle_\sigma$ so $\gamma=1$.
We conclude that 
$\langle \,,\,\rangle_\sigma$ is symmetric.
\hfill $\Box$ \\ 

\noindent At the end of Section 1
we indicated how finite-dimensional irreducible 
$\boxtimes$-modules give tridiagonal pairs.
The bilinear forms
in Theorem \ref{thm:bilex} are useful
in the theory of these tridiagonal pairs; we will
discuss this connection in a future paper.

\section{Suggestions for further research}

\noindent In this section we give some suggestions for further
research. We use the following notation.
For an integer $N\geq 1$ 
let 
$\mathfrak{sl}_2^{(N)}$ denote the Lie
algebra
$\mathfrak{sl}_2 \oplus 
\mathfrak{sl}_2  \oplus \cdots \oplus
\mathfrak{sl}_2$ 
($N$ copies).

\begin{problem}
\label{problem3}
\rm
For an integer $N\geq 1$ and
mutually 
distinct 
$a_1, a_2,\ldots, a_N$ in $\F\backslash \lbrace 0,1\rbrace$,
 consider the Lie algebra homomorphism
$\boxtimes \to 
\mathfrak{sl}_2^{(N)}$
that sends $\xi \mapsto (EV_{a_1}(\xi),EV_{a_2}(\xi),\ldots, EV_{a_N}(\xi))$
for $\xi \in \boxtimes$.
Describe the kernel of this homomorphism in terms of the
$\boxtimes$-generators
$\lbrace x_{ij} \,|\,i,j\in \I,\, i\not=j\rbrace$
and symmetric functions involving $a_1,a_2,\ldots, a_N$.
Also, find an attractive subset of $\boxtimes$
whose image under the homomorphism is a basis for
$\mathfrak{sl}_2^{(N)}$.
\end{problem}

\begin{example}
\rm
With reference to 
Problem  \ref{problem3},
assume $N=2$ and abbreviate $a=a_1$, $b=a_2$.
Then the kernel of the homomorphism
$\boxtimes \to 
\mathfrak{sl}_2
\oplus
\mathfrak{sl}_2
$ 
is generated as an ideal by
\begin{eqnarray*}
&&\lbrack x_{12},x_{03}\rbrack
-2ab(x_{23}-x_{01})-2(1-a)(1-b)(x_{31}+x_{02}),
\\
&&\lbrack x_{23},x_{01}\rbrack
-2(1-a^{-1})(1-b^{-1})(x_{31}-x_{02})
-2a^{-1}b^{-1}(x_{12}+x_{03}),
\\
&&\lbrack x_{31},x_{02}\rbrack
-2(1-a)^{-1}(1-b)^{-1}(x_{12}-x_{03})
-2ab(1-a)^{-1}(1-b)^{-1}(x_{23}+x_{01}).
\end{eqnarray*}
Also, the images
of
\begin{eqnarray*}
x_{12}, \quad
x_{23}, \quad
x_{31}, \quad
x_{01}, \quad
x_{02}, \quad
x_{03}
\end{eqnarray*}
under the homomorphism 
form a basis for
$\mathfrak{sl}_2
\oplus
\mathfrak{sl}_2$.
\end{example}

%\begin{problem}
%\rm
%\label{problem4}
%For an integer $N\geq 1$ find a generalization of Proposition
%\ref{lem:image}
%that applies to a surjective Lie
%algebra homomorphism
%$\boxtimes \to 
%\mathfrak{sl}^{(N)}_2$.
%\end{problem}

\begin{problem}
\label{problem1}
\rm Let $U(\boxtimes)$ denote the universal enveloping algebra
of $\boxtimes$. For $\sigma \in G$ find a formal sum
\begin{eqnarray*}
\Omega_\sigma = \sum_{n=0}^\infty t_n \qquad \qquad \qquad t_n \in 
U(\boxtimes)
\end{eqnarray*}
such that the following (i)--(iii) hold on each finite-dimensional
irreducible $\boxtimes$-module $V$:
\begin{enumerate}
\item
$t_n$ is 0 on $V$ for all but finitely many $n$.
\item $\Omega_\sigma$ is invertible on $V$.
\item $\sigma(\xi)-\Omega_\sigma \xi \Omega^{-1}_\sigma$ is 0 on
$V$ for all $\xi \in \boxtimes$.
\end{enumerate}
\end{problem}

\begin{note}
\rm Let $V$ denote a finite-dimensional irreducible
$\boxtimes$-module. Pick $\sigma \in G$ and consider
the sum $\Omega_\sigma$ from Problem
\ref{problem1}. Then the map
\begin{eqnarray*}
V &\mapsto & V
\\
v &\mapsto & \Omega_\sigma .v
\end{eqnarray*}
is an isomorphism of $\boxtimes$-modules from
$V$ to $V$ twisted via $\sigma$.
\end{note}

\begin{problem}
\label{problem2}
\rm
Find a short direct proof of Theorem
\ref{thm:modclass}.
\end{problem}

\begin{problem}
\label{prob:pd}
\rm
For an integer $N\geq 1$
let $V$ denote a vector space over $\F$ with
dimension $2N$.
%and let
%$\langle \,,\,\rangle $ denote a nondegenerate
%antisymmetric form on $V$.
%A subspace $U\subseteq V$ is called {\it isotropic}
%(with respect to 
%$\langle \,,\,\rangle$) whenever
%$\langle u,v\rangle=0$ for all $u,v \in U$.
Let $U_i$ $(i \in \I)$ denote 
$N$-dimensional subspaces of $V$
such that $U_i\cap U_j=0$ for
distinct $i,j\in \I$.
Show that there exists a $\boxtimes$-module structure
on $V$ such that 
$(x_{ij}+I)U_i=0$ and 
$(x_{ij}-I)U_j=0$
for distinct $i,j\in \I$.
\end{problem}

\begin{problem}
\rm
With reference to the $\boxtimes$-module $V$ in Problem
\ref{prob:pd}, show that the following are equivalent:
\begin{enumerate}
\item $V$ is a direct sum of irreducible $\boxtimes$-modules;
\item there exists a nondegenerate antisymmetric bilinear
form $\langle\,,\,\rangle $ on $V$ such that
 $\langle U_i,U_i\rangle=0 $ for $i \in \I$.
\end{enumerate}
\end{problem}

\begin{problem}
\label{prob:pdef}
\rm
With reference to Notation \ref{not},
assume $\F$ is the complex number field $\C$,
and that  $a \in \R$.
Further assume
that $\langle \eta_r,\eta_s\rangle \in \R$
for all distinct $r,s\in \I$.
Pick 
mutually distinct
$i,j,k,\ell\in \I$ and 
define $V_{\R} =\sum_{n=0}^d \R u_n$,
where $\lbrace u_n\rbrace_{n=0}^d$ is the basis
$\lbrack i,j,k,\ell\rbrack$ of
$V$ from
Lemma \ref{def:eta}.
 By
the data in Theorem
\ref{thm:trans} $V_{\R}$ is independent
of $i,j,k,\ell$. 
For a nonidentity $\sigma \in G$
consider the restriction of $\langle \,,\,\rangle_{\sigma}$
to 
$V_{\R}$. When is this restriction
positive definite?
\end{problem}

\begin{problem} \rm
Let $V$ denote a finite-dimensional irreducible
$\boxtimes$-module. Find all the linear transformations
$\psi:V\to V$ such that both
\begin{eqnarray*}
\lbrack x_{01},
\lbrack x_{01},
\lbrack x_{01},
\psi \rbrack \rbrack \rbrack= 
4 \lbrack x_{01},
\psi \rbrack,
\qquad \qquad 
\lbrack x_{23},
\lbrack x_{23},
\lbrack x_{23},
\psi \rbrack \rbrack \rbrack= 
4 \lbrack x_{23},
\psi \rbrack.
\end{eqnarray*}
\end{problem}

\begin{problem}
\rm
At the beginning of Section 7 we gave an action of $S_4$
on $\boxtimes$.   
Describe how
the $S_4$-module $\boxtimes$ 
decomposes into a direct sum of irreducible $S_4$-modules.
\end{problem}

\section{Acknowledgement}
\rm The authors thank Brian Curtin,
Eric Egge, Mark MacLean,
and Kazumasa Nomura for giving this paper a close reading
and offering many valuable suggestions.

\noindent Tatsuro Ito \hfil\break
\noindent Department of Computational Science \hfil\break
\noindent Faculty of Science \hfil\break
\noindent Kanazawa University \hfil\break
\noindent Kakuma-machi \hfil\break
\noindent Kanazawa 920-1192, Japan \hfil\break
\noindent email:  
{\tt tatsuro@kenroku.kanazawa-u.ac.jp} \hfil\break

\bigskip

\noindent Paul Terwilliger \hfil\break
\noindent Department of Mathematics \hfil\break
\noindent University of Wisconsin \hfil\break
\noindent 480 Lincoln Drive \hfil\break
\noindent Madison, WI 53706-1388 USA \hfil\break
\noindent email: {\tt terwilli@math.wisc.edu }\hfil\break

\end{document}